\documentclass[smallextended]{svjour3} 
\usepackage{amsfonts,amsmath,amstext,amssymb,amsopn}
\usepackage{graphicx}
\usepackage{color}
\usepackage{url}
\usepackage{algorithm}
\usepackage{algorithmic}
\usepackage{natbib}
\setcitestyle{numbers,comma,square}
\usepackage{tabularx}
\usepackage{tabulary,booktabs, threeparttable, stackengine,array}
\usepackage{mathtools}

\usepackage{authblk}
\usepackage[labelfont=bf]{subcaption}
\usepackage{bbold}

\begin{document}
\newcommand*{\affaddr}[1]{#1} 
\newcommand*{\affmark}[1][*]{\textsuperscript{#1}}
\newcommand{\exclude}[1]{}

\title{Optimization Models for Integrated Biorefinery Operations}
\author{Berkay Gulcan\affmark[1]     \and
		Sandra D. Eksioglu\affmark[2]   \and
		Yongjia Song\affmark[1]   \and
		Mohammad Roni\affmark[3]   \and
		Qiushi Chen\affmark[4] 
}

	
	
	\institute{B. Gulcan \at
		\email{bgulcan@clemson.edu}           
		\and
		\affaddr{\affmark[1]Department of Industrial Engineering, Clemson University, Clemson, SC, USA} \\
		\and
		\affaddr{\affmark[2]Department of Industrial Engineering, University of Arkansas, Fayetteville, AR, USA} \\
		\and
		\affaddr{\affmark[3]Idaho National Laboratory, Idaho Falls, ID USA}\\
		\and
		\affaddr{\affmark[4]Glenn Department of Civil Engineering, Clemson University, Clemson, SC, USA}\\}
	\date{Received: date / Accepted: date}
	\maketitle
\begin{abstract}
    Variations of physical and chemical characteristics of biomass lead to an uneven flow of biomass in a biorefinery, which reduces equipment utilization and increases operational costs. Uncertainty of biomass supply and high processing costs increase the risk of investing in the US's cellulosic biofuel industry. We propose a stochastic programming model to streamline processes within a biorefinery. A chance constraint models system's reliability requirement that the reactor is operating at a high utilization rate given uncertain biomass moisture content, particle size distribution, and equipment failure. The model identifies operating conditions of equipment and inventory level to maintain a continuous flow of biomass to the reactor. 
    The Sample Average Approximation method approximates the chance constraint and a bisection search-based heuristic solves this approximation. A case study is developed using real-life data collected at Idaho National Laboratory's biomass processing facility. 
    An extensive computational analysis indicates that sequencing of biomass bales based on moisture level, increasing storage capacity, and managing particle size distribution increase utilization of the reactor and reduce operational costs.\\
    \keywords{bioenergy \and biomass \and biorefinery operations optimization \and chance constraint optimization \and stochastic optimization \and system reliability}
\end{abstract}
\section{Introduction}\label{sec:Intro}

In the last two decades, significant investments have been made by a number of federal agencies to establish and develop US’s cellulosic biofuel industry. Despite of these investments, this industry remains a nascent concern. 
Major challenges faced by this industry are variations of biomass supply and variations of physical/chemical characteristics of biomass (such as, ash, moisture, carbohydrate contents, etc.), largely due to spatial variations in weather and soil, harvesting equipment used, etc. These variations lead to uneven flow of biomass in a biorefinery which affect equipment utilization rate and lead to inconsistent conversion rates. Variations in biomass supply and processing costs increase the risk of investing in this industry. 

A recent report from the US Department of Energy (DOE) identifies ``bulk solids handling and material flows through the system" as a critical component to achieve the design throughput of the conversion processes \citep{DOE2016}. The flow of materials is impacted by variations of biomass characteristics. For example, consider the scenario when a number of bales of different types of feedstock with different moisture level and ash content are processed on the same equipment (i.e., grinder). The resulting distribution of particle size and particle uniformity of processed biomass varies from one bale to the next. These variations negatively affect the integration of biomass feeding system and conversion process, which lead to  low/unreliable on-stream time of equipment and low utilization of the reactor.

The main objective of this work is to improve the reliability of biomass feeding system in a biorefinery via a stochastic optimization model. The reliability of the biomass feeding system is defined as the probability of achieving the targeted reactor utilization rate under stochastic biomass characteristics and stochastic equipment failure. In this study, we only consider equipment failures caused by  clogging, overflow (in conveyors), and overheat. These failures are mainly due to  biomass characteristics, such as, moisture content and particle size distribution. We focus on the reactor because it is the most expensive equipment in a biorefinery.  Achieving the targeted utilization rate at all time is expensive and not practical, because it requires the use of biomass with low moisture and ash continents. Therefore, it is economical to focus on achieving the utilization of reactor at a target rate, most of the time. We model this requirement using chance constraints.  

In order to optimize the performance of the reactor, the proposed optimization model identifies optimal equipment operating conditions, inventory level for a given sequence of bales based on biomass moisture levels. These approaches align with strategies presented in the literature, which indicate that, system's reliability can be improved by increasing redundancy, sequencing components and products,  etc. \cite{DavidCoit_reliability_opt2018}. In our system, bale sequencing impacts the size of inventory, and inventory holding mitigates disruption of biomass flow due to  failure and operating conditions of equipment. 

 Biomass density is one of the main factor that impacts the biomass flow in the system. Studies show that biomass particle size distribution and moisture content significantly affect  biomass density~\cite{ZhouCornStover2008},~\cite{Crawford_2016_Flowability_CornStover}. To evaluate these relations we develop a Discrete Element Method (DEM) model. DEM is a computational model which describes the mechanical bulk behaviour of granular materials \cite{DEM1979}. DEM simulates the movement and interaction of particles with each other and with the system. The results of this simulation are used to develop regression functions which calculate the density of biomass as a function of moisture level and particle size distribution. We use biomass density to calculate biomass flow in the system. The proposed optimization uses a multi-period network flow model to capture the flow of biomass in the system, and the flow of biomass into the reactor during the planning horizon.   

We develop a case study using historical data from the Process Development Unit (PDU), a biomass facility, at Idaho National Laboratory (INL). We use this case study to validate  the proposed model and conduct numerical experiments. In particular, our numerical experiments are designed to evaluate the impact of  ($i$) biomass characteristics (e.g., moisture level, particle size distribution and density) on achieving the targeted reactor utilization rate and minimizing operating costs; ($ii$) equipment failure on reactor utilization rate and costs; and ($iii$) strategies, such as, sequencing bales, changing operating condition of equipment, and keeping inventory, on reducing the risk of achieving the targeted reactor utilization rate.

We expect that the results of this study will help biorefineries to develop strategies  which lead to increased reactor utilization and minimize costs; and to identify optimal operational condition in face of stochastic biomass characteristics and equipment failure. These outcomes will facilitate  commercial scale generation of biofuels at competitive cost. In the long run, these outcomes will strengthen the sustainable bioeconomy of US, enhance the security of energy supplies, reduce dependencies on fossil fuels, and reduce greenhouse gas (GHG) emissions. A strong sustainable bioeconomy has  additional socioeconomic benefits, such as generating new green jobs, growth of rural economy and social stability, among others~\cite{DOMAC2005Bioenergy,Sims2003Bioenergy,You2012Biofuel}.

The remainder of the paper is organized as follows. In Section~\ref{sec:Review}, we review the literature. In Section~\ref{sec:Problem}, we describe the  
problem and the modeling approach. A case study is presented in Section~\ref{sec:Case}. In Section~\ref{sec:Results}, we present numerical results and analysis. We close with concluding remarks in Section~\ref{sec:Conclusion}.

\section{Literature Review} \label{sec:Review}
The two main streams of literature closely related to this study are system reliability and optimization of biorefinery operations. There are a number of studies in the literature that focus on biomass supply chain optimization~\cite{Chen_Onal_2014,EKSIOGLU20091342, Memisoglu2015, RONI2014115}, which impact the availability of biomass in a biorefinery. In this paper we do not consider supply chain decisions and assume that the mix and quantity of biomass to process is given as an input to the problem. 
\subsection{System Reliability} \label{sec:SystemReliability}
Maximizing system's reliability is a relevant problem in manufacturing and service systems which operate under uncertainty. During the last decade, a number of researchers studied the problem and developed solution approaches. The field of system reliability optimization is still expanding. A recent  survey paper by \cite{DavidCoit_reliability_opt2018} presents a comprehensive review of this literature. It provides a chronological grouping of the research in the field into three eras, the era of mathematical programming, era of pragmatism, and era of active reliability improvement. It also classifies the problems studies in this field into four major groups ($i$) redundancy allocation, ($ii$) reliability allocation, ($iii$) reliability-redundancy allocation, and ($iv$) assignment and sequencing. Each group of research provides strategies to improve the system’s reliability. The model we propose contributes to \emph{redundancy allocation} and \emph{sequencing} streams of research since it identifies the amount of inventory of biomass  and the sequencing of bales processed which minimize the cost of meeting a targeted system reliability level.      

A wide range of methods are used to solve reliability problems. For example, \cite{Bellman34, Fyffe1986Reliability,Ghare1969Reliability} use dynamic programming to solve a redundancy allocation problem. Work by \cite{Misra_IP-Reliability_1991, Prasad_ReliabilityOpt_2000} use linear programming and integer programming models to minimize costs or  maximize reliability. Many system reliability problems are combinatorial in nature due to component selection and sequencing. Additionally, models of complex systems use nonlinear constraints and objectives. Therefore, exact solution approaches are computationally expensive. To address these computational challenges, a number of metaheuristics are proposed, such as, genetic algorithms and simulated annealing, which are shown to be efficient and robust. For example, \cite{Painton1995PCReliability} developed a genetic algorithm to model a personal computer component configuration problem,  \cite{YANG1999Reliability_GA} use a genetic algorithm to design a nuclear power plant. While metaheuristics do not guarantee to find the global optima, researchers showed that metaheuristics are efficient and robust methods to solve complex system reliability problems.

Recent developments capture the impact of uncertainty to the reliability of a system. For example, \cite{Basciftci2018,Li_ChanceConstr_Reliability_2008} use chance-constrained stochastic programming to model and solve a system reliability problem under uncertainty. Similar to this work, our proposed model uses  chance constraints  to capture the impact of biomass characteristics and equipment failure on the utilization of the reactor.  
 
Current efforts in this field are focused on the development of data driven optimization models \cite{YANG2010671, Alsina2018, MENG2020773}. The proposed work also contributes to this stream of research by using sensor-based data collected from the PDU to develop the DEM models and estimate failure probabilities of equipment.

\subsection{Optimization of Biorefinery Operations}\label{sec:optBio}
Biorefinery operations include biomass pre-processing, biomass handling and biomass storage. Most of the literature related to biorefinery operations focuses on the evaluation of  design parameters of the equipment used \cite{Crawford_2016_Flowability_CornStover,DAI2012716}. Other studies analyze energy consumption of an equipment as a function of its design parameters and biomass characteristics \cite{osti_1369631, osti_1133890, yancey2015size}. These studies are limited in scope since they do not capture the interactions among equipment and the impact of equipment on the performance of the system as a whole. 

A number of studies use optimal control models to minimize the energy consumption of an equipment. For example, \cite{NUMBI20161653} propose a deterministic, non-linear optimization model that identifies optimal operating parameters of a crusher. By introducing a time-to-use electricity tariff, the model achieves additional cost savings. Work by \cite{ZHANG20101929} presents an optimal control model for a series of conveyor belts. The authors compare two different operational structures, one in which the conveyors' speed is fixed, and another in which the conveyors' speed is adjusted over time. The comparison of these models points to the benefits of having the flexibility to update the operating conditions of equipment. To summarize, optimal control models are complex non-linear programs, thus, are used to model only parts of a system. Using such an approach to model the performance of a system under uncertainty, would lead to a notoriously complex model. This is the main reason why we did not follow such an approach.  

A few studies present models for an optimal design of operations in a biorefinery. For example, \cite{Viet2012BiorefOpt} propose a two steps approach to optimize the design and operations of a biorefinery.  This model identifies an optimal biorefinery configuration for a given set of biomass feedstocks and available conversion technologies. The first step, called the bi-directional synthesis, identifies the intermediary chemicals that can be produced using the feedstocks available. The second step optimizes a network flow problem that identifies the mix of intermediary chemicals and corresponding conversion technologies that minimize the total production cost. Work by \cite{ZONDERVAN2011BiorefOpt} uses a mixed-integer and nonlinear network optimization model to identify an optimal design for a biorefinery. The model identifies process sequences to optimize  production of a set of biofuels and bioproducts. Both works,  \cite{Viet2012BiorefOpt} and \cite{ZONDERVAN2011BiorefOpt}, use deterministic models since they consider fixed biomass characteristics, flow rates and yields.  Different from these works, our research models uncertainties in this system. The use of our proposed models leads to robust decisions.  

A number of studies use DEM models to study the flowability of materials. 
DEM simulates movements and interactions of particles with each other and with the system.  By using the law of force, DEM calculates the rotation and velocity of each particle and simulates the movement of material in an equipment (or system) realistically \cite{Hohner2012}. 
Additionally, DEM models consider the impact of bond strengths of particles and deformability of materials to flowability \cite{XIA20191}.  The main reason for the growing interest in these models is the lack and inaccuracy of experimental data necessary to model the flowability of granular material \cite{XIA20191}. In the last decade DEM models have been used to evaluate the flowability of different biomass feedstocks \cite{XIA20191, guo2020discrete}. A study by \cite{SCHERER2016896} develops a DEM model to simulate the drying process of wood chips in a rotary dryer. They analyzed the impacts of equipment design on the process rate and quality. \cite{OREFICE2017347} develops a DEM model of a horizontal screw conveyor, a common transportation equipment for biorefineries. They study the impacts of screw rotational velocity and initial filling level of the conveyor on the flow rate of the material. \cite{XIA20191} developed a DEM model to study the impacts of deformable pinewood chips in cyclic loading test. This study \cite{XIA20191} shows that DEM can be used to simulate particle deformation, which is critical to the performance of equipment used in a biorefinery (e.g., grinders and pelleting equipment). Work by \cite{guo2020discrete} present a bonded-sphere DEM model designed for switchgrass particles. This model generates data which are used to evaluate  the relationship between biomass density and moisture level and particle size distribution. Our proposed work uses the results of this model to estimate material density. This integration of the results of the DEM model within the optimization model is one of the contributions of the proposed work.

\section{Problem Description and Formulation}\label{sec:Problem} The goal of the proposed mathematical model is to identify process controls and inventory levels which minimize system wide costs while meeting the desired biomass specifications for biochemical conversion; and achieving the designed production throughput. Our model is developed using data from the PDU since it provides a good representation of the processes used for feeding of biomass to the reactor in a biorefinery. Next, we describe the problem, present a mathematical formulation, present a solution approach and model extensions.    

\subsection{Problem Description} 
\begin{figure}[H]
	\centering
	\includegraphics[width=0.8\linewidth]{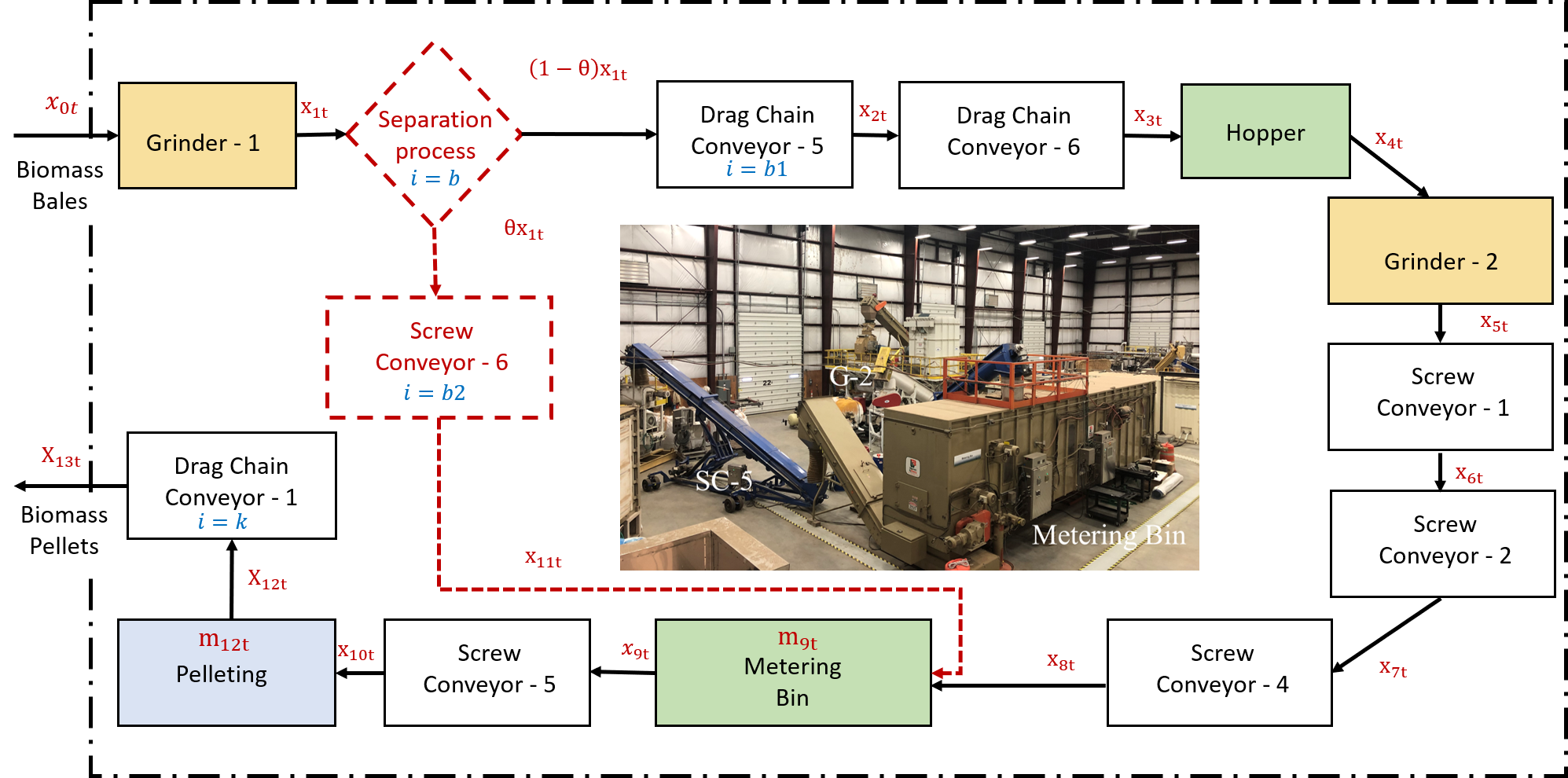}
	\caption{Biomass Feeding Processes of PDU}
	\label{fig:PDU}
\end{figure}

Figure \ref{fig:PDU} presents the flowchart of the different processes of PDU. Bales of biomass enter the system, one at a time, via a conveyor belt. Biomass is processed in Grinder 1. Processes biomass undergoes a process which separates biomass based on particle size. Large biomass particles are transported via conveyors to Grinder 2 for further processing. Small biomass particles are transported to the Metering bin for storage. Biomass  processed in Grinder 2 also is stored in the Metering bin. From the Metering bin biomass is transported to the pelleting machine. Pelleted biomass is fed to the reactor.


We group the equipment in this system based on their tasks into processing, transportation and storage. A processing equipment converts biomass from its original format (e.g. baled, log, or coarse-shredded) to its final format (e.g., ground biomass and pellet). These equipment are grinders 1 and 2 and the pelleting mill. Most equipment of PDU are transportation equipment, such as, conveyors. Storage equipment, such as, the metering bin, store in-process inventory.

The following are the  \emph{assumptions} we make, which help us modeling this system. First, historical data from PDU reports only the average moisture level of bales processed. This is calculated using data collected via sensors located in the conveyor belt which feeds biomass bales to the system. It takes a few minutes to process one bale. Since the moisture level of a bale is heterogeneous, we generate the moisture level of biomass in different time periods using a uniform distribution. The mean of this distribution equals the average moisture level of the corresponding bale. The lower and upper bounds of the distribution are established via discussions with experts. 

Second, the moisture level of biomass is reduced after being processed in an equipment. Historical data is used to estimate these reductions in moisture level for each processing equipment.

Third, work by \cite{yancey2015size} shows that the probability of clogging increases with moisture level of biomass being processed. Thus, in practice, the processing speed of equipment is adjusted based on the moisture level of biomass. In our model we establish upper bounds on the amount of biomass that can be processed in a equipment based on the moisture level of biomass.

\subsection{Problem Formulation} 
We use stochastic optimization to model the uncertainties observed in the system. This model is a two-stage stochastic program with chance constraints. We use a two-stage stochastic program because some of the decisions are made before uncertainties (such as, biomass moisture content, particle size, and bulk density) reveal. These first-stage decisions include the  processing speed of equipment over the entire planning horizon, $\mathbf{V}$, and  the initial inventory level in the metering bin,  $\mathbf{I}_0$. Next, given a realization of the biomass characteristics $\omega$ and the first-stage decisions, the second-stage decisions are made, which include biomass flows in the system, $\mathbf{X}(\omega)$, and the inventory level in the metering bin, $\mathbf{I}(\omega)$. 
We use a chance constraint to model the reliability level of the reactor. Since maintaining the reactor operating at the targeted utilization rate all the time is expensive, we aim to achieve this target utilization rate most of the time. Let $r$ represent this target rate and $1-\epsilon$ the desired reliability of the reactor. We use a chance constraint to ensure that the reactor achieves the target rate $r$ at least $(1-\epsilon)*100$ percent of the time.   
  
We use network flow constraints to model the flow of biomass within the system in every period during the planning horizon $\mathcal{T}$. Let $\bf G=(\bf N, \bf A)$, with a node set $\bf N$ and arc set $\bf A$, represent the network structure of the system. The set of nodes represent equipment of the PDU (Figure \ref{fig:PDU}), and the set of arcs represent the flow of biomass from one equipment to next. The following is a list of notation we use.

\begin{table}[h]
	\scriptsize
		\begin{tabular}{ll}
			\toprule
			\multicolumn{2}{c}{\textbf{SETS:}}\\
		    ${\bf M}$ & The set of moisture levels of biomass.\\
		    $\mathcal{T}$ & The set of time periods in the planning horizon, $\mathcal{T} := \{1,2,\ldots, T\}$\\
		    ${\bf E}^p$ & The set of processing equipment.\\
		    ${\bf E}^r$ & The set of transportation equipment.\\
		    ${\bf E}^m$ & The set of storage equipment.\\
		    ${\bf N} =  {\bf E}^p \cup {\bf E}^r \cup {\bf E}^m$\\
			\midrule
			\multicolumn{2}{c}{\textbf{PARAMETERS:}}\\
			$A$ & Node arc incidence matrix of $\bf G$.\\
			$A^{\prime}$ & Node arc incidence matrix of $\bf G^{\prime}(\bf {\bf E}^r, \bf A^{\prime})$.\\
			$\mathbf{h}$ & Inventory holding cost for storage equipment (in \$/ton).\\
	        $f(\cdot)$ & Energy consumption and operational cost function.\\
	        $\kappa := \{\kappa_t\}_{t\in \mathcal{T}}$ & Moisture level of the biomass bales \\
			$\kappa_t$ & Moisture level of the biomass bales in time period t ($\kappa_t \in M$).\\	
			$\bar{v}(\kappa_t)$ & Upper bound of equipment processing speed.\\ 
			$\bar{\iota}(\kappa_0)$ & Initial inventory holding capacity in the storage equipment.\\
			$r$ & The target reactor feeding rate (in tons).\\
			$\epsilon$ & Risk tolerance parameter.\\
			\midrule
			\multicolumn{2}{c}{\textbf{RANDOM PARAMETERS:}}\\
			$\bar{\iota}(\omega)$ & The inventory holding capacity in the storage equipment for processing biomass $\omega$.\\
			$\underline{\iota}(\omega)$ & The inventory lower threshold in the storage equipment for processing biomass $\omega$.\\
			$\omega$ & A random vector that contains the stochastic biomass characteristics.\\
			$\Omega$ & The support of the probability distribution of $\omega$ \\
			\midrule
			\multicolumn{2}{c}{\textbf{DECISION VARIABLES:}}\\
			$\mathbf{V} = \{\mathbf{V}_t\}_{t \in \mathcal {T}}$,  & Equipment processing speed (in meters or rotations per time period).\\ 
			where $\mathbf{V}_t := \{V_{it}\}_{i\in {\bf N}}, \ \forall t \in \mathcal {T}$ &
			\\
			$\mathbf{I}_0 := \{I_{i0}\}_{i\in {\bf E}^m}$ & Initial inventory level in the storage equipment \\
			$\mathbf{I}(\omega) := \{\mathbf{I}_t(\omega)\}_{t \in \mathcal {T}}$ & Inventory level when processing biomass with characteristics $\omega$ (in  tons).\\
			$\mathbf{I}_t(\omega) := \{I_{it}(\omega)\}_{i\in {\bf E}^m}$ & Inventory level when processing biomass with characteristics $\omega$ at period $t$ (in  tons).\\
			$\mathbf{X}(\omega) = \{\mathbf{X}_t(\omega)\}_{t \in \mathcal {T}}$ & Biomass flow when processing biomass with characteristics $\omega$ (in tons).\\
			$\mathbf{X}_t(\omega) := \{X_{it}(\omega)\}_{i\in {\bf N}}$ & Biomass flow when processing biomass with characteristics $\omega$ at period $t$ (in tons).\\
			\bottomrule
			\end{tabular}
\end{table}%
For simplicity of presentation, we next present a succinct formulation $(P)$ using generic functional notation such as $f(\cdot)$, $b(\cdot)$ and $d(\cdot)$ to model the relationship between decision variables and random variables. A detailed formulation is provided in Appendix B. Given a sequence of bales with moisture level $\kappa$ to process in the planning horizon of $\mathcal{T}$, the proposed stochastic program  is given by:

\begin{alignat}{2}
(P) \hspace{0.5in} \min: \ & 
\mathcal{Z} = \mathbf{h}^\top\mathbf{I}_0 + \mathbb{E}\left [f(\mathbf{X}(\omega),\mathbf{I}(\omega), \omega) \right] \label{obj}\\
\text{s.t. } &   A^{\prime}\mathbf{X}(\omega) = 0, \ \forall \omega\in \Omega, \label{eqn:Flow1}\\
& 0 \leq \mathbf{V}_t \leq \bar{v}(\kappa_t), \ \forall t \in \mathcal{T},
\label{eqn:Speed_bound}\\
& \mathbf{X}_t(\omega) = b(\mathbf{V}_t,\mathbf{I}_0,\mathbf{I}_t(\omega) , \omega), \ \forall t \in \mathcal{T}, \omega \in \Omega, \label{eqn:Flow2}\\ 
&  \mathbf{X}_t(\omega) \le d(\mathbf{V}_t, \omega), \ \forall t \in \mathcal{T}, \omega \in \Omega, \label{eqn:Flow3}\\
& \underline{\iota}(\omega) \leq \mathbf{I}_t(\omega) \leq \bar{\iota}(\omega), \ \forall t \in \mathcal{T}, \omega \in \Omega, \label{eqn:Inv_bound}\\
& 0 \leq \mathbf{I}_0 \leq \bar{\iota}(\kappa_0), \label{eqn:Inv_bound_init}\\
 &  \mathbb{P} \left( R(\mathbf{X}(\omega), \omega) \geq r \right) \geq 1-\epsilon.  \label{eqn:Probabilistic}
\end{alignat}

The objective function~\eqref{obj} computes the total expected cost of the system, which includes the initial inventory holding cost and the expected operating costs. The  operating costs, $f(\cdot)$, include the cost of energy consumed by equipment during the planning horizon. 

Constraints \eqref{eqn:Flow1} represent the flow balance constraints of transportation equipment. Constraints \eqref{eqn:Speed_bound} set an upper bound $\bar{v}(\kappa_t)$ on the processing speed of each equipment. These bounds depend on the moisture level $\kappa_t$ of  biomass being processed at period $t$. 

Constraints \eqref{eqn:Flow2} calculate the flow from storage and processing equipment. For storage equipment, these represent the inventory balance constraints. Pelleting mill, different from other processing equipment, has a positive residence time. Therefore, an inventory balance constraint is used to calculate the flow of biomass from this equipment. The inventory balance constraint is written explicitly as follows:
\begin{align*}
    & I_{i,t-1}(\omega) + \sum_{j \in \boldsymbol{\delta_i^-}}X_{jts}(\omega) = I_{it}(\omega) + X_{it}(\omega), \ \forall i\in {\bf E}^m, t\in\mathcal{T}\setminus\{1\},\\
    & I_{i,0} + \sum_{j \in \boldsymbol{\delta_i^-}}X_{jts}(\omega) = I_{i1}(\omega) + X_{i1}(\omega), \ \forall i\in {\bf E}^m,
\end{align*}
where $\boldsymbol{\delta}^-_i$ represents the set of equipment that feeds into equipment $i$.  

Constraints \eqref{eqn:Flow2} are simpler for other processing equipment, such as grinders.  
Let $\tilde{d}_{it}(\omega)$ represent the density of biomass after processed in grinder $i$, and let $\gamma_i$ represent the cross section area of the discharge opening of grinder $i$, thus $\gamma_i \tilde{d}_{it}(\omega)$ is the mass discharging rate of the grinder. Given the infeed rate of the grinder $V_{it}$, the flow from a grinder is calculated as:  
\[ X_{it}(\omega) = \gamma_i \tilde{d}_{it}(\omega) V_{it}.\]

Constraints \eqref{eqn:Flow3} represent the upper limit on the amount of biomass flow from  equipment $i\in {\bf N}$. For example, let $i$ be a screw conveyor with cross section area equal to $\gamma_i$. Thus, $\gamma_i \tilde{d}_{it}(\omega)$ is the mass discharging rate, and  $\gamma_i \tilde{d}_{it}(\omega) V_{it}$ is the maximum amount of biomass that can be discharged by the conveyor. 

\[ X_{it}(\omega) \leq \gamma_i \tilde{d}_{it}(\omega) V_{it}.\]

Constraints \eqref{eqn:Inv_bound} and \eqref{eqn:Inv_bound_init} set upper and lower bounds in the amount of biomass that is stored in the metering bin. 
These bounds are  random in~\eqref{eqn:Inv_bound} because the mass of biomass that can be stored depends on its density, which is a random parameter. We set lower bounds since this is a practice used at PDU to maintain consistent flow of biomass from the metering bin. 
Finally,  \eqref{eqn:Probabilistic} is the chance constraint. 
 
\subsection{A Sample Average Approximation of (P)}\label{sec:SAA_model}
Chance-constrained stochastic programs are typically difficult to solve. This is because, for a given solution, it is often difficult to check its feasibility since it may require multi-dimensional integration. Furthermore, the feasible region defined by chance constraints is in general non-convex. To address these challenges, many researchers use the Sample Average Approximation (SAA) method to approximate these constraints. The corresponding  model is  an integer program which is easier to solve using commercial solvers~\citep{atlason2008optimizing, Luedtke_Ahmed_2008, Pagnoncelli_Ahmed_Shapiro_2009}.

The SAA model replaces the true probability distribution of random parameters with an empirical distribution obtained from random samples. Let $\{\omega_1, \omega_2, \cdots, \omega_S\}$ be a set of $S$ independent and identically distributed realizations (scenario) of $\omega$, which are obtained via a Monte Carlo simulation. Thus, the probability associated with each of these realizations equals $\frac{1}{S}$. In addition, let the second-stage decision variables $\mathbf{X}_t(\omega_s)$ and $\mathbf{I}_t(\omega_s)$ be denoted by short-handed notation $\mathbf{X}_{ts}$ and $\mathbf{I}_{ts}$, respectively, representing the flow value and inventory level in each period $t$ under each scenario $s$. Let $\mathbb{1}(\cdot)$ be an indicator function. The SAA approximation of (P) can be written as:  
\begin{align}
(\hat{P}) \hspace{0.5in} \mathcal{Z}^S = \min \ & \mathbf{h}^\top \mathbf{I}_0 + \frac{1}{S}\sum_{s=1}^{S} f(\mathbf{X}_{ts},\mathbf{I}_{ts}, \omega_s) \nonumber\\
\text{s.t. } &  \eqref{eqn:Speed_bound}, \eqref{eqn:Inv_bound_init}, \nonumber \\
 &   A^{\prime}\mathbf{X}_{s} = 0, \  s= 1, 2,\ldots, S, \label{eqn:sFlow1}\\
 & \mathbf{X}_{ts} = b(\mathbf{V}_t,\mathbf{I}_0,\mathbf{I}_{ts}, \omega_s), \ \forall t \in \mathcal{T},  s= 1, 2,\ldots, S, \label{eqn:sFlow2}\\ 
 &  \mathbf{X}_{ts} \le d(\mathbf{V}_t, \omega_s), \ \forall t \in \mathcal{T}, s= 1, 2,\ldots, S, \label{eqn:sFlow3}\\
 & 0 \leq \mathbf{I}_{ts} \leq \bar{\iota}(\omega_s), \ \forall t \in \mathcal{T}, s= 1, 2,\ldots, S, \label{eqn:sInv_bound}\\
 & \frac{1}{S}\sum_{s=1}^{S}\mathbb{1}\left[R(\mathbf{X}_s, \omega_s) \geq r \right]\geq 1-\hat{\epsilon}.  \label{eqn:sProbabilistic}
\end{align}

In formulation $(\hat{P})$, the value of $\hat{\epsilon}$ may be different from the true risk parameter ${\epsilon}$ used in formulation ($P$). 
Based on ~\cite{Luedtke_Ahmed_2008}, when $\hat{\epsilon} < \epsilon$, the probability that a feasible solution of $(\hat{P})$ is feasible to $({P})$ increases with the sample size $S$. In our experiment we choose to use $\hat{\epsilon} < \epsilon$ in experiments where we consider equipment failures  and use $\hat{\epsilon} = \epsilon$ elsewhere.

Constraint \eqref{eqn:sProbabilistic} of ($\hat{P}$) is not linear because it includes an indicator function. One could resort to an integer programming reformulations, however, it can be time-consuming to solve \cite{Luedtke_Ahmed_2008}. Instead, we use a heuristic approach for solving ($\hat{P}$) which penalizes the violation of the chance  constraint \cite{charnes1955optimal, abdelaziz2012solution}. To this end, we introduce variables $\boldsymbol{\mathcal{U}}$ to quantify the amount of violations, and a penalty parameter $\pi > 0$. The following is a linear approximation of ($\hat{P}$).  
\begin{align}{}
(\bar{P}) \hspace{0.5in} \bar{\mathcal{Z}}^S = \min \ &  
\mathbf{h}^\top \mathbf{I}_0 + \frac{1}{S}\sum_{s=1}^{S} f(\mathbf{X}_s,\mathbf{I}_s, \omega_s)
+ \sum_{s=1}^{S}\pi \mathcal{U}_s  \nonumber\\
\text{s.t. } &  \eqref{eqn:Speed_bound}, \eqref{eqn:Inv_bound_init}, \eqref{eqn:sFlow1} - \eqref{eqn:sInv_bound}, \nonumber \\
 & R(\mathbf{X}_s, \omega_s)
 + \mathcal{U}_s - \mathcal{J}_s  = r, \ s = 1,2,\ldots, S, \label{eqn:ssProbabilistic}\\
& \mathcal{U}_s, \mathcal{J}_s  \geq 0,  \ s = 1,2,\ldots, S.
\end{align}
The new term in the objective function of $(\bar{P})$ penalizes the difference between the average amount of biomass fed to the reactor and the target value $r$, if the target cannot be achieved.

For a given value of the penalty parameter $\pi > 0$, formulation $(\bar{P})$, which is  a two-stage stochastic linear program, is easy to solve. Note that, parameter $\pi$ is not known in advance and one needs to identify an appropriate $\pi$ value so that the resulting solution satisfies the chance constraint $(1-\hat\epsilon)*100$ of the time. It is obvious that $\pi \rightarrow \infty$  leads to solutions for which $\mathcal{U}_s=0$. In this case, the chance constraint is satisfied for all $s$. Conversely, $\pi \rightarrow 0$ leads to solutions for which $\mathcal{U}_s\geq 0$. In this case, the chance constraint may be satisfied in less then $[1 - \hat{\epsilon}]S$ of the scenarios generated. Thus, we design a bisection search algorithm which identifies the smallest value of $\pi$ for which  the number of scenarios which satisfy the chance constraint (i.e. $R(\mathbf{X}_s, \omega_s) \geq r$) is close to $[1 - \hat{\epsilon}]S$. We present this algorithm in detail in the Appendix D. 

We use the stochastic Benders decomposition approach to solve the corresponding two-stage stochastic linear program $(\bar{P})$ \cite{Birge_Louveaux_1997}. We use Benders decomposition because the problem size increases with the number of scenarios. We use the multi-cut version of Benders decomposition.
 
\subsection{Model Extensions}\label{sec:Extensions}
{\bf Integration of the DEM Models.} One of the DEM models developed by \cite{guo2020discrete} simulates  grinders 1 and 2. This model is calibrated and validated using historical data from the PDU. The model estimates the bulk density of biomass via simulation. The simulation results are used to develop regression functions that capture the relationship between  biomass density (dependent variable), and moisture level  and particle size distribution (independent variables). Let $\rho^{j}_i$ represent the $j$-th percentile of the particle size distribution of biomass processed in equipment $i$. The regression model is given by:
\begin{align*}
\displaystyle \tilde{d}_{it} = \alpha^0_i + \alpha^1_i m_{it}+ \alpha^2_i \rho^{50}_i + \alpha^3_i \frac{\rho^{90}_i}{\rho^{10}_i} + \tilde{\xi}_i,
\end{align*}
where $\tilde{\xi}_i$ corresponds to a random error term following a normal distribution with mean equal to zero and a constant standard deviation. This regression model is used to compute the density of biomass in the proposed model for each of the scenarios generated.   
 
{\bf Incorporating equipment failure.} The proposed model can also be extended by incorporating random equipment failures. We consider two types of equipment failures: short-duration and long-duration. Short-duration failures are  due to the overfeeding of the system, overflowing of conveyors, or overheating of the grinders. Long-duration failures are due to clogging of an equipment, which typically happens when processing biomass with high moisture level and large particle size. We assume that short-duration and long-duration failures are independent of each other. 

For a short-duration failure, we assume that the system restarts automatically within a few seconds. A long-duration failure does typically require an operator to unclog the equipment and restart the system. These failures are observed in equipment located in the upstream of the metering bin (the storage equipment), such as the grinders and the corresponding conveyors. Failures are not observed in the downstream of the metering bin since the flow of biomass is controlled via the inventory. Additionally, the moisture level of biomass is reduced considerably by the time it leaves the metering bin, and particle size of biomass is reduced considerably after being processed in grinder 2. Both factors reduce the probability of clogging.    

To model the random equipment failures, we assume that the time between consecutive equipment failures and time to repair are random variables following certain probability distributions. We assume that the time between failures follows a Weibull distribution and the duration of a failure follows a uniform distribution. The parameters used for these distributions are summarized in Table~\ref{tab:stoppage}, which are validated using data collected from the PDU. We generate an offline operational schedule of equipment (with up and down time) over the planning horizon using random samples according to these distributions. These random operating schedules are incorporated in the model via parameters $o_{it}(\omega)\in \{0,1\}$, where, $o_{it}(\omega) = 0$ if  equipment $i$ is down in period $t$, and $o_{it}(\omega) = 1$ if the equipment is operating. We update constraints~\eqref{eqn:Flow2} which calculate the flow of biomass from an equipment. 
\[ X_{it}(\omega) = o_{it}(\omega)\gamma_i \tilde{d}_{it}(\omega) V_{it}.\]

\section{Case Study}\label{sec:Case}
\subsection{Data Collection} 
We develop a case study using  historical data and the design of the biomass processing system from the PDU. The operational costs are collected from the Biomass Logistics Model (BLM), which is an integrated software framework that simulates the entire supply chain and calculates associated costs, energy consumption and GHG emissions  \cite{Cafferty2013BLM}. The summary table of the operational costs is presented in Appendix A.

Our data set summarizes the data related to the processing of 14 bales of switchgrass during a period of 4 days at PDU that includes the amount of biomass processed, equipment throughput, electrical current and power consumption per time period (of 0.2 seconds). The moisture level of biomass during each time period is recorded via sensors. Changes of moisture level and dry matter losses during grinding and pelleting operations, and system stoppages are also recorded. 

We use another data set from the PDU that includes the density of biomass and particle size distribution measured at three different points in this process: ($i$) before the process begins while biomass is still in bale format, ($ii$) after biomass is processed in the first grinder, and ($iii$) after  biomass is processed in the second grinder. The PDU uses a screen of size 76.2mm (3 inch) in the primary grinder and screens of size 6.35mm ($\frac{1}{4}$ inch) in the separation process and the secondary grinder in order to separate particles based on their length. This data is considered in modeling of the separation process.

Work by~\cite{HANSEN2019BiomassSC}  provides additional data about bulk density of switchgrass when different harvesting equipment are used. They report  densities which vary from $171.26$ to $234.81$ $kg/m^3$. 

Finally, we use data generated by the DEM model to create a regression model that represents the relationship between moisture content, particle size distribution and bulk density. Tables~\ref{tab:Cost} to \ref{tab:infeed_limit} in the Appendix summarize the input data used. We have consulted the experts and operators of the PDU during model development, verification and validation.

\subsection{Data Analysis}
We assume that the density of biomass bales follows a triangular distribution. This distribution is typically used when the number of samples is small and conducting additional sampling is expensive~\cite{BeyondBeta, Thomopoulos2017}, which is the case here.

The moisture level of bales in our case study varies from $5\%$ to $25\%$, and bales are grouped into low (5\% to 10\%), medium (10\% to 17.5\%) and high (17.5\% to 25\%) moisture levels. We fitted a uniform distribution to describe the distribution of moisture in each level. 

The following regression models are developed using data from the DEM model. The DEM model assumes that biomass particle size follows a uniform distribution. Regression \eqref{reg1} presents biomass density after processed at grinder 1. Regression \eqref{reg2} presents biomass density after processed at grinder 2. These regressions estimate density as a function of moisture level and particle size. The error term $\xi_{it}$ is normally distributed, $\boldsymbol{\xi^1_{it}$  $\sim$ $N(0,3.106)}$ and $\boldsymbol{\xi^2_{it}$  $\sim$ $N(0,10.783)}$. 

\begin{equation}\label{reg1}
\tilde{d}_{it} = 56.183 + 65.312 \tilde{m}_{it}  - 8.473\rho^{50}_i + 0.015\frac{\rho^{90}_i}{\rho^{10}_i} + \tilde{\xi}^1_{it},\\
\end{equation}
\begin{equation}\label{reg2}
\tilde{d}_{it} = 186.348 + 206.1697 \tilde{m}_{it} - 110.302 \rho^{50}_i + 0.709\frac{\rho^{90}_i}{\rho^{10}_i} + \tilde{\xi}^2_{it}.
\end{equation}

Table \ref{tab:reg_stats} in Appendix E summarizes the statistical analysis of regression functions \eqref{reg1} and \eqref{reg2}. The $R^2$ values for these regression functions are over 94\%. Based on the P-values found, the impact of $\frac{\rho^{90}_i}{\rho^{10}_i}$ on biomass density is statistically insignificant.  Thus, we do not include $\frac{\rho^{90}_i}{\rho^{10}_i}$ in regression functions \eqref{reg1} and \eqref{reg2} used in our numerical analysis. 

The time-to-failure and duration of an equipment failure depend on the moisture level of biomass. We use historical data from the PDU to estimate these parameters. The duration of a failure is modeled using the uniform distribution, and time-to-failure is modeled using the Weibull distribution. The corresponding parameters for these distributions are summarized in Table \ref{tab:stoppage}.

\begin{table}[htp!]
	\scriptsize
	\centering
	\begin{minipage}{0.45\textwidth}
	{%
		\begin{tabular}{lcccc}
			\toprule
			
			& \multicolumn{2}{c}{\textbf{Time-to-failure}}& \multicolumn{2}{c}{\textbf{Failure Duration}}\\
			& \multicolumn{2}{c}{\bf Weibull}& \multicolumn{2}{c}{\bf Uniform}\\
			\cmidrule(lr){2-3}
			\cmidrule(lr){4-5}
			\textbf{Moisture} & \textbf{Shape}& \textbf{Scale} & \textbf{Min.}& \textbf{Max.}\\
			\textbf{ \hspace{.1mm} Level}& &$(sec)$ &$(sec)$ &$(sec)$\\
			\cmidrule(lr){1-1}
			\cmidrule(lr){2-2}
			\cmidrule(lr){3-3}
			\cmidrule(lr){4-4}
			\cmidrule(lr){5-5}
			Low		&  $1.16$ & $9.94$ & $0.0$ & $4.0$  \\
			Medium	&  $0.83$ & $15.09$ & $0.0$ & $7.0$  \\
			High	&  $0.59$ & $22.91$ & $0.0$ & $12.3$  \\			
			\bottomrule
			\\
			\multicolumn{5}{c}{\textbf{Short-duration Equipment Failure}}
	\end{tabular}
	}
	\end{minipage}\hfill
	\begin{minipage}{0.45\textwidth}
	{
	\begin{tabular}{lcccc}
			\toprule
			
			& \multicolumn{2}{c}{\textbf{Time-to-failure}}& \multicolumn{2}{c}{\textbf{Failure Duration}}\\
			& \multicolumn{2}{c}{\bf Weibull}& \multicolumn{2}{c}{\bf Uniform}\\
			\cmidrule(lr){2-3}
			\cmidrule(lr){4-5}
			\textbf{Moisture} & \textbf{Shape}& \textbf{Scale} & \textbf{Min.}& \textbf{Max.}\\
			\textbf{ \hspace{.1mm} Level}& &$(sec)$ &$(sec)$ &$(sec)$\\
			\cmidrule(lr){1-1}
			\cmidrule(lr){2-2}
			\cmidrule(lr){3-3}
			\cmidrule(lr){4-4}
			\cmidrule(lr){5-5}
			Low		&  $5.50$ & $95$ & $20.0$ & $35.0$  \\
			Medium	&  $5.00$ & $90$ & $20.0$ & $45.0$  \\
			High	&  $4.50$ & $85$ & $20.0$ & $60.0$  \\			
			\bottomrule
			\\
			\multicolumn{5}{c}{\textbf{Long-duration Equipment Failure}}
	\end{tabular}
	}
	\end{minipage}
	\caption{Distribution Parameters}
	\label{tab:stoppage}
\end{table}%

\section{Numerical Experiments and Sensitivity Analysis}\label{sec:Results}
The goal of our numerical experiments is two fold. First, in Section  \ref{sec:model-evaluation} we evaluate the performance of the model and algorithms proposed. Next, in Sections~\ref{sec:base-analyis} to \ref{sec:particleSize} we evaluate the performance of the system.   The performance of the system is measured via ($i$)  reactor utilization, ($ii$)  operational cost, and ($iii$)  inventory level. Section~\ref{sec:base-analyis} summarizes the results of our base-case scenario and Sections \ref{sec:BaleProc} to \ref{sec:particleSize}  summarize the sensitivity analysis with respect to bale sequencing, equipment failure rate, storage capacity, and biomass characteristics, respectively. 

\subsection{Evaluating the Performance of the Model and Algorithms Proposed}\label{sec:model-evaluation}
In our numerical experiments, we present solutions of SAA model ($\hat{P}$) under different settings. 
We solve ($\hat{P}$) using a bisection search-based heuristic approach (see Algorithm \ref{alg:BinearySearch} in Appendix D). To find an appropriate sample size $S$ for the SAA model ($\hat{P}$), we conducted the following stability test. We first solved model ($\hat{P}$) by varying the number of scenarios. Each scenario represents a sample path realization of random biomass characteristics ($\omega$) over the planning horizon. For a given number of scenarios, we ran $10$ replications and computed the relative difference between the objective values associated with these replications. Our experiments showed that with 100 scenarios, the relative difference between the maximum and minimum objective values among the $10$ replications is only $0.3\%$, indicating that the sample size of $100$ is appropriate. We thereby use a sample of size $ S = 100$ scenarios in our numerical experiments.

To evaluate the performance of the solutions found via the SAA model $(\hat{P})$, we consider an out-of-sample evaluation procedure for the first-stage solution $\mathbf{V}^*$ and $\mathbf{I}_0^*$. 
We independently generate $10,000$ scenarios in the out-of-sample test and check if solution ($\mathbf{V}^*$ and $\mathbf{I}_0^*$) satisfies the chance constraint~\eqref{eqn:Probabilistic} with these $10,000$ scenarios. Specifically, we compute the flow values $\mathbf{X}_s$ according to equation~\eqref{eqn:sFlow2} using the first-stage solution $\mathbf{V}^*$ and $\mathbf{I}_0^*$ for each out-of-sample scenario $s$, and then we check if the reactor utilization target is achieved for at least $(1-\epsilon)\times 10000$ scenarios. 

Notice that, constraints~\eqref{eqn:Flow1}-\eqref{eqn:Inv_bound_init} of $(P)$ define system capacities and thresholds that need to be satisfied under all possible realizations of the random variables. Although in the SAA model, the corresponding constraints~\eqref{eqn:sFlow1} to~\eqref{eqn:sInv_bound} are satisfied for the scenarios generated, enforcing these constraints only for the scenarios used by the SAA model may not be sufficient to ensure that the obtained first-stage solution $\mathbf{V}^*$ and $\mathbf{I}_0^*$ is feasible under any possible scenario. To ensure that the first-stage solution to the SAA model is ``absolutely'' feasible under any scenario, one may need to resort to a robust optimization approach by modeling the entire support of random variables $\omega$ as an uncertainty set. This will result in additional (deterministic) constraints in the SAA model ($\hat{P}$), however, we do not expect that this will affect the analysis below.

We next demonstrate the value of the proposed chance-constrained stochastic programming model. First, we compare the solution of $(\hat{P})$ with that of the \emph{mean-value (MV) problem} based on their out-of-sample performance. The MV problem considers a single deterministic scenario where all the random parameters in the problem are replaced by their mean value for the entire duration of the planning horizon. Our experiment result shows that the MV solution ($\mathbf{V}$ and $\mathbf{I}_0$) led to infeasible solutions in all 10,000 scenarios of the out-of-sample evaluation, and the amount of violation in the metering bin capacity constraint was as much as 666\% in some of the scenarios. This indicates that it is necessary to use a stochastic programming model like the proposed model ($P$) and its SAA model ($\hat{P}$) to adequately address the stochasticity in the problem. 

We next justify the use of the proposed chance-constrained stochastic program to model system's reliability, which is defined as the probability of achieving the desired target feeding rate (see constraint~\eqref{eqn:Probabilistic}). Based on our discussions with the PDU operator, we set the reliability level to be 90\%. We conducted a sensitivity analysis to evaluate the impact of the established  reliability level on system's performance. The results of this analysis are summarized in Table~\ref{tbl:SystemReliability} and Figure~\ref{fig:FlowvsRisk}. These experiments assume that the biomass processed in the system has low moisture level. 
Similar observations are made when the system processes biomass with higher moisture level. 

Table~\ref{tbl:SystemReliability} and Figure~\ref{fig:FlowvsRisk} illustrates the trade-off between costs, target reactor feeding rate, average inventory, maximum inventory and $\epsilon$. As reliability level decreases, the operational cost, the average inventory and the maximum inventory decrease and a higher reactor target rate is achieved.

\begin{minipage}{\textwidth}
  \begin{minipage}[b]{0.49\textwidth}
  \captionsetup{justification=centering}
    \centering
    \scriptsize{
    \begin{tabular}{cccc}
			\toprule
			\textbf{System}&\textbf{Reactor}&\textbf{Average}&\textbf{Maximum} \\
			\textbf{Reliability}&\textbf{Target}&\textbf{Inventory}&\textbf{Inventory} \\
		    $(\%)$ & $(dt/ hr)$ & $(tons)$ & $(tons)$\\
			\cmidrule(lr){1-1}
			\cmidrule(lr){2-2}
			\cmidrule(lr){3-3}
			\cmidrule(lr){4-4}
			100 & 3.68 & 0.56 & 0.88 \\
			99 & 3.76 & 0.56 & 0.87\\
			90 & 3.79 & 0.54 & 0.84\\
			80 & 3.81 & 0.52 & 0.81\\
			\bottomrule
	\end{tabular}}
      \captionof{table}{System's Performance vs.  Reliability}
	\label{tbl:SystemReliability}
  \end{minipage} 
  \hfill
  \begin{minipage}[b]{0.49\textwidth}
    \centering
    \includegraphics[width=0.9\linewidth]{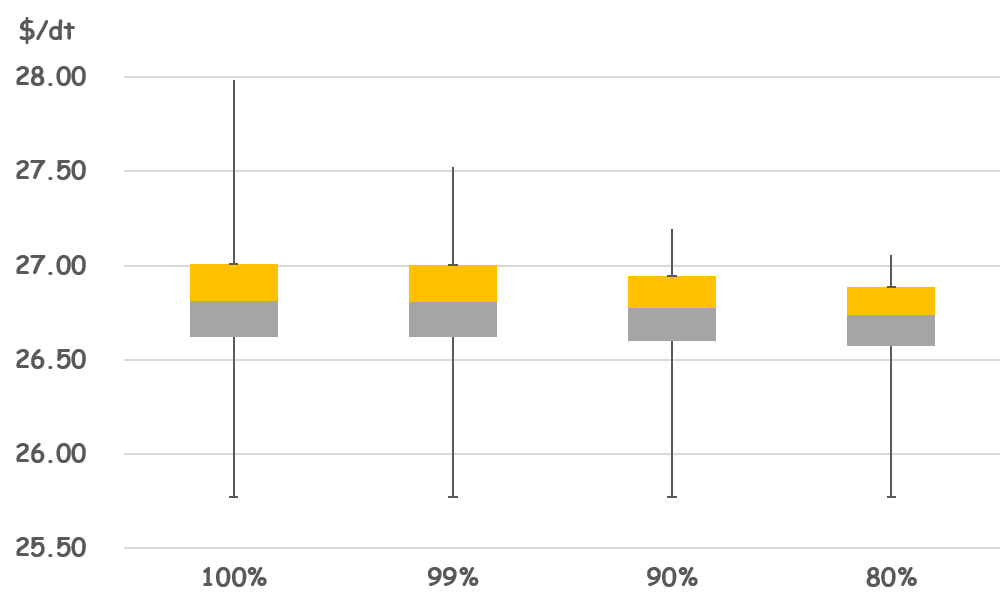}
    \captionof{figure}{Total Cost vs. Reliability}
    \label{fig:FlowvsRisk}
  \end{minipage}
\end{minipage}

\subsection{Base Case Analysis}\label{sec:base-analyis}
In this section, we summarize our experimental results for a base-case problem, which considers that every biomass bales, each of the same moisture level, is processed in a system with no equipment failures. We consider that  storage capacity is $49.1\ m^3$, the particle size follows a uniform distribution, and the risk level is set to 90\%. Additional problem parameters are presented in Appendix C. 


Tables~\ref{tbl:Cost_base} and \ref{tbl:Inv_base} summarize the performance of the system under two moisture levels and two reactor capacities. Throughout this section, we consider 2.7 dt/hr to be  low capacity, and 4.8 dt/hr to be high. We observe that reactor utilization is higher and the operational cost is lower when processing low moisture biomass since the in-feed rate is higher.  The inefficiency of the system when operating high moisture biomass lead to low reactor utilization. For example, when reactor's capacity is high, processing high moisture biomass leads to only 44\% utilization of the reactor.   

\begin{table}[htbp]
\centering
	\scriptsize
	{%
		\begin{tabular}{ccccccc}
			\toprule
			\textbf{Bale Moisture}&\textbf{Reactor} &\textbf{Reactor}&\textbf{Reactor}& \textbf{Energy} &\textbf{Fixed}&\textbf{Total} \\
			\textbf{Level}&\textbf{Capacity} &\textbf{Flow}&\textbf{Utilization}&\textbf{Cost} &\textbf{Costs}&\textbf{Cost}\\
		    & $(dt/ hr)$ & $(dt/ hr)$ & $(\%)$ & $(\$/ dt)$ & $(\$/ dt)$ & $(\$/ dt)$\\
			\cmidrule(lr){1-1}
			\cmidrule(lr){2-2}
			\cmidrule(lr){3-3}
			\cmidrule(lr){4-4}
			\cmidrule(lr){5-5}
			\cmidrule(lr){6-6}
			\cmidrule(lr){7-7}
			Low & Low & 2.25 & 83 & 3.85 & 46.60 & 50.45\\
            Low & High & 3.85 & 80 & 2.25 &	27.26 &	29.54 \\
            High &	Low  & 2.10&	77	&6.96	&50.33	&57.29\\
            High &	High  & 2.11&	44 &6.93	&50.13	&57.06\\
			\bottomrule
	\end{tabular}}
	\caption{Base-Case Problem: Reactor Utilization \& Costs}
	\label{tbl:Cost_base}
\end{table}%

The results of Table~\ref{tbl:Inv_base} indicate that, when processing biomass with high moisture level, the initial inventory level is higher and infeed rate is lower then when processing biomass with low moisture level. The additional inventory is needed to ensure a continuous flow of biomass to the reactor when moisture level is high. The corresponding average and maximum inventory are also higher.


\begin{table}[htbp]
\centering
	\scriptsize
	{%
		\begin{tabular}{cccccccc}
			\toprule
			\textbf{Bale Moisture}&\textbf{Reactor}&\textbf{Reactor} &\textbf{Initial}&\textbf{Average}& \textbf{Maximum} &\textbf{Average System}&\textbf{Metering Bin} \\
			\textbf{Level}&\textbf{Capacity}&\textbf{Flow} &\textbf{Inventory}&\textbf{Inventory}&\textbf{Inventory} &\textbf{Infeed Rate}&\textbf{Conveyor Speed}\\
		    & $(dt/ hr)$ & $(dt/ hr)$ & $(tons)$ & $(tons)$ & $(tons)$ & $(inch/min)$ & $(inch/min)$\\
			\cmidrule(lr){1-1}
			\cmidrule(lr){2-2}
			\cmidrule(lr){3-3}
			\cmidrule(lr){4-4}
			\cmidrule(lr){5-5}
			\cmidrule(lr){6-6}
			\cmidrule(lr){7-7}
			\cmidrule(lr){8-8}
			Low & Low &2.25	&0.23	&0.37	&0.53	&6.7&4.0\\
            Low & High &3.85	&0.25	&0.56	&0.88	&11.5&6.9\\
            High &	Low & 2.10	&2.12	&0.71	&2.12	&5.5 & 4.4\\
            High &	High & 2.11	&2.18	&0.54	&2.18	&5.5 & 4.5\\
			\bottomrule
	\end{tabular}}
	\caption{Base-Case Problem: Inventory Level \& Equipment Setting}
	\label{tbl:Inv_base}
\end{table}%

The results of Tables~\ref{tbl:Cost_base} and \ref{tbl:Inv_base} provide the best and worst performance of the system since we consider that biomass moisture level is either low or high. Section \ref{sec:BaleProc} evaluates system's performance when a mix of bales of different moisture level are processed. 

\subsection{Sensitivity Analysis: Bale Sequencing
}\label{sec:BaleProc}
In this section, we analyze the impact of bale sequencing on the performance of the system. We consider that $60\%$ of the bales processed are of low moisture, $10\%$ are of medium moisture, and $30\%$ are of high moisture. We sequence bales based on moisture level. The following sequences are considered, long, short and random. The long sequence considers that every bale of a particular moisture level is processed before the processing of bales of another moisture level. In our experiments, the long sequence begins by processing  high moisture bales, then medium moisture bales and ends with low moisture bales. Such a sequence leads to the worst performance among all long sequences because it requires higher initial inventory to maintain a continuous feeding of the reactor while processing high moisture bales. We consider this specific sequence to highlight the difference in the performance of long versus short sequences. 

%

The short sequence follows a pattern of $60\%$ low,  $10\%$ medium, and $30\%$ high moisture bales. This pattern repeats itself several times during the planning horizon. Finally, a random sequence  processes bales of different moisture level using a random pattern of high, medium and low moisture. The results presented in Table \ref{tbl:Cost_sequences}  are the averages from 10 problem generated using these sequencing approaches.

\begin{table}[htbp]
\centering
	\scriptsize
	{%
		\begin{tabular}{ccccccc|cc}
			\toprule
			\textbf{Sequencing}&\textbf{Reactor} &\textbf{Reactor}&\textbf{Reactor}& \textbf{Energy} &\textbf{Fixed}&\textbf{Total} & \textbf{Average}  & \textbf{Max}\\
			\textbf{Approach}&\textbf{Capacity} &\textbf{Flow}&\textbf{Utilization}&\textbf{Cost} &\textbf{Costs}&\textbf{Cost}&\textbf{Inventory}&\textbf{Inventory}\\
		     & $(dt/ hr)$ & $(dt/ hr)$ & $(\%)$ & $(\$/ dt)$ & $(\$/ dt)$ & $(\$/ dt)$& $(tons)$& $(tons)$\\
			\cmidrule(lr){1-1}
			\cmidrule(lr){2-2}
			\cmidrule(lr){3-3}
			\cmidrule(lr){4-4}
			\cmidrule(lr){5-5}
			\cmidrule(lr){6-6}
			\cmidrule(lr){7-7}
			\cmidrule(lr){8-8}
			\cmidrule(lr){9-9}
			Long & Low &  2.31	&85	 &4.60 	 &45.60 &50.20 &0.58 &1.71\\
            Long & High & 3.39	&71	 &3.13 	 &31.14 &34.26 &0.72 &2.07\\
            Short &	Low  & 2.30	&85	 &4.60 	 &45.62 &50.23 &0.41 &0.61\\
            Short &	High  &3.68 &77	 &2.88 	 &28.61 &31.49 &1.05 &1.76\\
            Random &Low  & 2.31	&85	 &4.60 	&45.61 	&50.21 	&0.41 &0.63\\
            Random &	High  &3.62	&75	 &2.93 &29.17 &32.11 &0.84 &1.92\\
			\bottomrule
	\end{tabular}}
	\caption{Bale Sequencing: Results of the Sensitivity Analysis
	}
	\label{tbl:Cost_sequences}
\end{table}%
 
Table \ref{tbl:Cost_sequences} summarizes the results of the sensitivity analysis with respect to bale sequencing. Short sequences outperform  long and random sequences. Short sequences perform best when the capacity of the reactor is high. In this case, the flow to the reactor is highest at 3.68 dt/hr, utilization is highest at 77\%, and total cost is lowest at \$31.49/dt. 

 
\begin{figure}[htbp]
\centering
\begin{subfigure}{.5\textwidth}
  \centering
  \includegraphics[width=.9\linewidth]{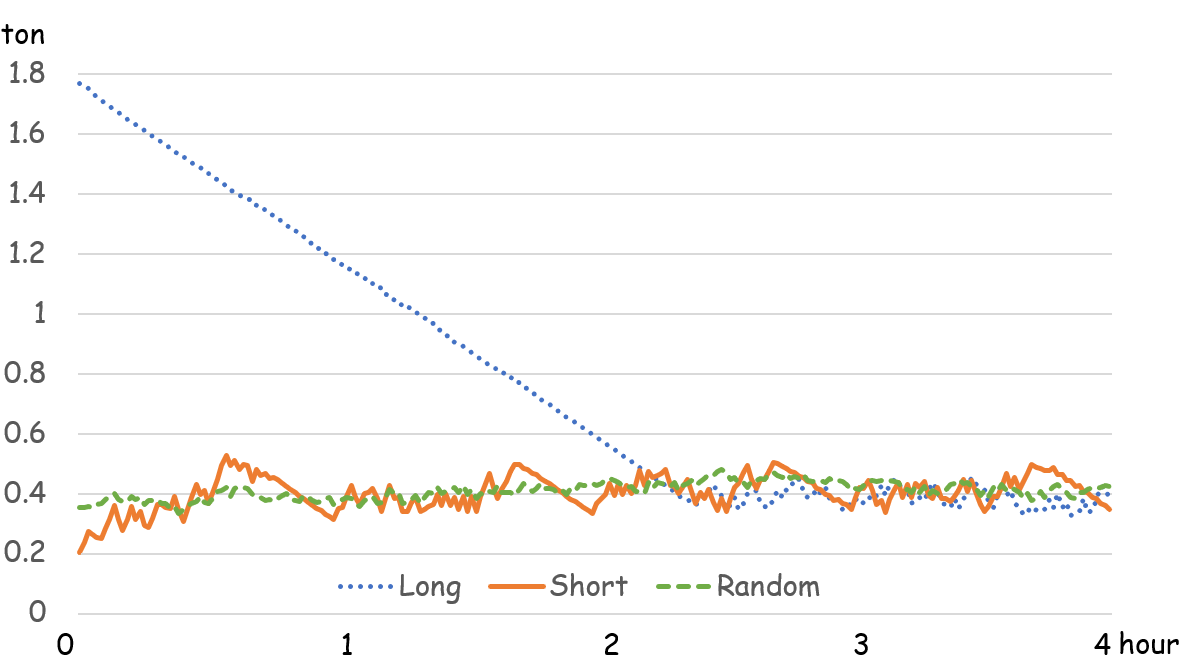}
  \caption{Reactor Capacity: Low}
  \label{fig:Inv-sub1}
\end{subfigure}%
\begin{subfigure}{.5\textwidth}
  \centering
  \includegraphics[width=1\linewidth]{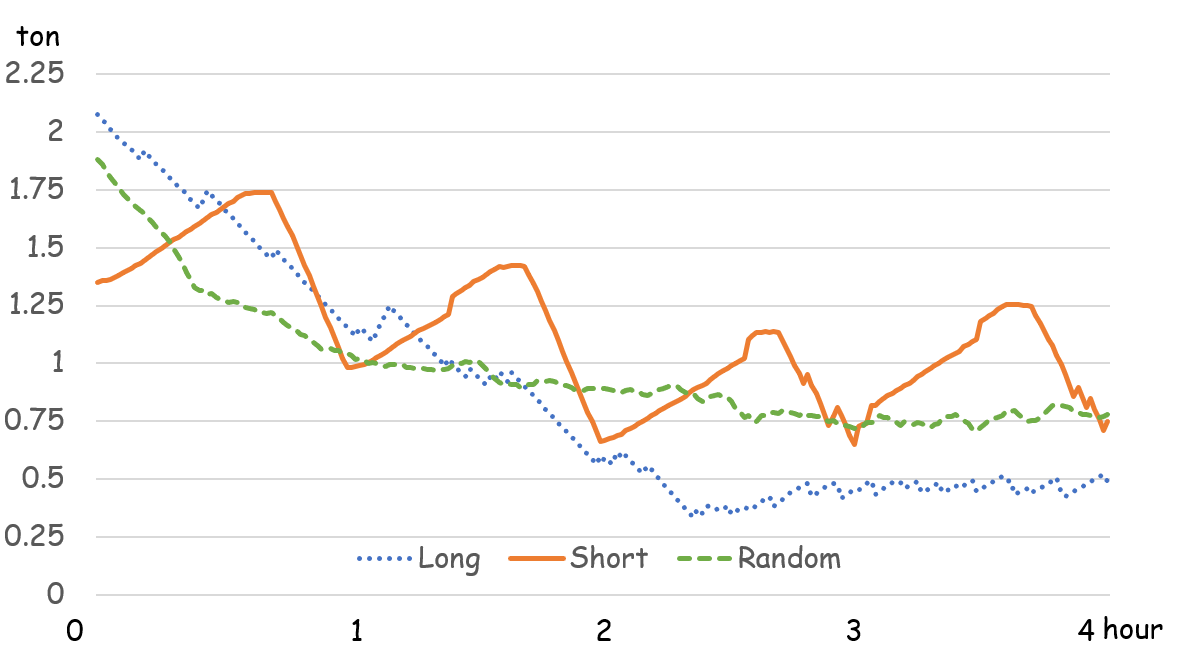}
  \caption{Reactor Capacity: High}
  \label{fig:Inv-sub2}
\end{subfigure}
\caption{Inventory Level for Different Sequences}
\label{fig:Inv-seq}
\end{figure}

 Figure~ \ref{fig:Inv-seq} presents the inventory level of each sequence. 
In both, Figures~\ref{fig:Inv-sub1} and Figure~\ref{fig:Inv-sub2}, one can
observe a cyclic pattern of inventory level for short sequences since inventory accumulates at a high rate while processing low moisture bales, inventory accumulates at a low rate while processing medium moisture bales, and inventory is used while processing high moisture bales. Long sequences begin with the highest level of initial inventory in order to maintain a continuous flow of biomass to the reactor since these sequences begin with high moisture bales. The performance of the random sequence is similar to the short sequence when the processing capacity of the reactor is low. The short sequence outperforms the random sequence when the processing capacity of the reactor is high. The short sequence leads to about 2\% lower total costs. Nevertheless, our model does not consider the cost of sequencing the bales. Thus, if the costs of creating a short sequence is high, a biomass processing plant should ensure that bales are processed at some random patterns. The plant should avoid the use of long sequences.         


	
\subsection{Sensitivity Analysis: Equipment Failure}
In this section, we analyze the impact of short-duration and long duration  equipment failure on the performance of the system.   

\noindent {\bf Short-duration Failures:} Recall that, short-duration  failures are due to the overfeeding of the system, overflowing of conveyors, or overheating of the grinders. These failures last for no more than 15 seconds, after which the equipment begins working automatically. Table~\ref{tbl:Cost_shortfail} summarizes the impact of short-duration failures on reactor utilization, total cost and inventory level. 

A comparison of the results of Table~\ref{tbl:Cost_shortfail} with Tables~\ref{tbl:Cost_base} and \ref{tbl:Inv_base}, indicates that short duration failures do not impact reactor's utilization and total cost when every bale has low moisture. However, the average inventory level increases by 13.6\% and maximum inventory level increases by 15.9\%. These increase of inventory enables the system to maintain a continuous flow of biomass to the reactor.  The performance of the system deteriorates most when every bale processed has high moisture level. The decrease in reactor utilization is 4\% and the increase in costs is 4.2\%.  


\begin{table}[htbp]
\centering
	\scriptsize
	{%
		\begin{tabular}{cc|cc|ccc}
			\toprule
					\textbf{Biomass Feeding}&\textbf{Reactor} &\textbf{Reactor} &\textbf{Total} &\textbf{Initial} &\textbf{Average} &\textbf{Maximum} \\
			\textbf{Pattern}&\textbf{Capacity} &\textbf{Utilization}&\textbf{Cost} &\textbf{Inventory}&\textbf{Inventory} &\textbf{Inventory}\\
		    &  &$(\%)$ & $(\$/dt)$&  (tons) & (tons) & (tons)\\
			\cmidrule(lr){1-1}
			\cmidrule(lr){2-2}
			\cmidrule(lr){3-3}
			\cmidrule(lr){4-4}
			\cmidrule(lr){5-5}
			\cmidrule(lr){6-6}
			\cmidrule(lr){7-7}
             Low	 & Low & 83 & 50.49	&0.23	&0.42	&0.62\\
             Low	 & High &80 &29.54	&0.71	&0.61	&0.94\\
             High &	Low & 74 & 59.69	&2.15	&1.10	&2.15 \\
             High &	High & 42 & 59.41	&2.15	&0.64	&2.15 \\
           \cmidrule(lr){1-1}
			\cmidrule(lr){2-2}
			\cmidrule(lr){3-3}
			\cmidrule(lr){4-4}
			\cmidrule(lr){5-5}
			\cmidrule(lr){6-6}
			\cmidrule(lr){7-7}
            Long &	Low & 84 & 50.59	&2.20	&0.68	&2.20\\
            Long &	High &70 &34.72	&2.20	&0.72	&2.20\\
            \cmidrule(lr){1-1}
			\cmidrule(lr){2-2}
			\cmidrule(lr){3-3}
			\cmidrule(lr){4-4}
			\cmidrule(lr){5-5}
			\cmidrule(lr){6-6}
			\cmidrule(lr){7-7}
            Short &	Low & 84 & 50.74	&0.35	&0.48	&0.76\\
            Short &	High & 75 & 32.38	&2.37	&0.98	&2.37\\
            \cmidrule(lr){1-1}
			\cmidrule(lr){2-2}
			\cmidrule(lr){3-3}
			\cmidrule(lr){4-4}
			\cmidrule(lr){5-5}
			\cmidrule(lr){6-6}
			\cmidrule(lr){7-7}
            Random & Low & 85& 50.22	&0.40	&0.50	&0.79\\
            Random & High& 74 & 32.73	&2.10	&0.78	&2.10\\
 			\bottomrule
	\end{tabular}}
	\caption{System Performance under Short-duration Equipment Stoppages}
	\label{tbl:Cost_shortfail}
\end{table}%


\noindent{\bf Long-duration Failures:} Recall that, long-duration failures are due to clogging which typically happens when biomass with high moisture level and large particle size is processed. Table~\ref{tbl:Cost_longfail}  summarizes the impact of long-duration failures on reactor utilization, total cost and  inventory level. The frequency and duration of clogging depends on the moisture level of biomass being processed (Figure \ref{tab:stoppage}). Thus, every problem in this table faces a different frequency and duration of long failures.

A comparison of the results of Table~\ref{tbl:Cost_longfail}  with Table~ \ref{tbl:Cost_sequences} indicate that long duration failures greatly impact the performance of the system. The utilization of the reactor is reduced by 5 to 58\%. The total cost is increased 5 to 140\%. The maximum inventory level reaches its limits in all but the problem where every bale has low moisture level. These observations raise the question whether an increase of storage capacity would allow for additional accumulation of the inventory which could be processed by the reactor during the time an equipment is down. Section \ref{sec:CapacityIncrease} investigates the impact of increasing storage capacity on system's performance. We also observe that the performance of the system is worst when every bale has high moisture level. The performance of the system is better when short sequences are processed as compared to long or random sequences.

\begin{table}[htbp]
\centering
	\scriptsize
	{%
		\begin{tabular}{cc|cc|ccc}
			\toprule
					\textbf{Biomass Feeding}&\textbf{Reactor} &\textbf{Reactor} &\textbf{Total} &\textbf{Initial} &\textbf{Average} &\textbf{Maximum} \\
			\textbf{Pattern}&\textbf{Capacity} &\textbf{Utilization}&\textbf{Cost} &\textbf{Inventory}&\textbf{Inventory} &\textbf{Inventory}\\
		    &  &$(\%)$ & $(\$/dt)$&  (tons) & (tons) & (tons)\\
			\cmidrule(lr){1-1}
			\cmidrule(lr){2-2}
			\cmidrule(lr){3-3}
			\cmidrule(lr){4-4}
			\cmidrule(lr){5-5}
			\cmidrule(lr){6-6}
			\cmidrule(lr){7-7}
			Low	    &Low & 79 & 53.00	&0.29	&1.19	&2.03\\
            Low	    &High &52 & 45.68	&1.26	&1.24	&2.12\\
            High	&Low &44 &99.48	&2.15	&1.41	&2.19\\
            High	&High & 25 &99.37	&2.22	&1.33	&2.22\\
           \cmidrule(lr){1-1}
			\cmidrule(lr){2-2}
			\cmidrule(lr){3-3}
			\cmidrule(lr){4-4}
			\cmidrule(lr){5-5}
			\cmidrule(lr){6-6}
			\cmidrule(lr){7-7}
            Long	&Low	&61 & 69.84	&2.15	&1.40	&2.20\\
            Long	&High   & 38 & 63.68 &2.15	&1.43	&2.24 \\
            \cmidrule(lr){1-1}
			\cmidrule(lr){2-2}
			\cmidrule(lr){3-3}
			\cmidrule(lr){4-4}
			\cmidrule(lr){5-5}
			\cmidrule(lr){6-6}
			\cmidrule(lr){7-7}
            Short	&Low	&63 & 67.43	&1.44	&1.25	&2.03\\
            Short	&High   &40& 59.77	&1.63	&1.24	&2.07\\
            \cmidrule(lr){1-1}
			\cmidrule(lr){2-2}
			\cmidrule(lr){3-3}
			\cmidrule(lr){4-4}
			\cmidrule(lr){5-5}
			\cmidrule(lr){6-6}
			\cmidrule(lr){7-7}
			Random	&Low	&56& 75.95	&1.72	&1.32	&2.10\\
            Random	&High   &31& 76.81	&1.65	&1.06	&2.01\\	 
			\bottomrule
	\end{tabular}}
	\caption{System Performance under Long-duration Equipment Failures}
	\label{tbl:Cost_longfail}
\end{table}%

\subsection{Sensitivity Analysis: Storage Capacity} \label{sec:CapacityIncrease}
In this section, we analyze the impact of an increase of storage capacity on the performance of the system. We consider a $25\%$ and $50\%$ increase of capacity, and resolve  ($i$) the base-case problem, ($ii$) problems with short-duration equipment failures, and ($iii$) problems with long-duration equipment failures. Next, we summarize our findings. 

\noindent {\bf Base-case Problem:}  Table~\ref{tbl:Capacity-NoStop3} summarizes the impact of increased storage capacity on costs and reactor utilization for the cases when  reactor's capacity is low and high. The last two columns of this table present the change in these performance measures as compared to the results of the base-case problem (Tables~\ref{tbl:Cost_base} and \ref{tbl:Cost_sequences}).

\begin{table}[htbp]
\centering
	\scriptsize
	{%
		\begin{tabular}{cc|cc|cc}
			\toprule
			\textbf{Biomass Feeding} &\textbf{Inventory Capacity}&\textbf{Reactor}& \textbf{Total}&\textbf{Reactor}& \textbf{Total}\\
			\textbf{Pattern} &\textbf{Increase}&\textbf{Utilization}&\textbf{Cost}&\textbf{Utilization $\Delta$}&\textbf{Cost $\Delta$}\\
		     & $(\%)$ & $(\%)$ & $(\$/dt)$& $(\%)$ & $(\%)$\\
			\hline\\
			\multicolumn{6}{c}{\bf Low Reactor Capacity}\\
			\hline
             High  & 25 & 79 & 56.18& 2.1 & -1.9\\
             High  & 50 & 81& 54.69& 5.0 & -4.5\\
             \hline\\
             	\multicolumn{6}{c}{\bf High Reactor Capacity}\\
            \hline
			 High		&25	&45 & 55.71 &2.5	&-2.4\\
             High		&50	&46 &54.60  &4.7	&-4.3\\
             \cmidrule(lr){1-1}
			\cmidrule(lr){2-2}
			\cmidrule(lr){3-3}
			\cmidrule(lr){4-4}
			\cmidrule(lr){5-5}
			\cmidrule(lr){6-6}
			Long		&25	&72 &33.75            &1.7	&-1.5\\
             Long		&50	&73 &33.30           &3.1	&-2.8\\
            \cmidrule(lr){1-1}
			\cmidrule(lr){2-2}
			\cmidrule(lr){3-3}
			\cmidrule(lr){4-4}
			\cmidrule(lr){5-5}
			\cmidrule(lr){6-6}
            Short		&25	     &77 & 31.15           &1.1	&-1.1\\
            Short		&50	     &77 &31.15           &1.1	&-1.1\\
            \cmidrule(lr){1-1}
			\cmidrule(lr){2-2}
			\cmidrule(lr){3-3}
			\cmidrule(lr){4-4}
			\cmidrule(lr){5-5}
			\cmidrule(lr){6-6}
			Random		&25	     &76 &31.68          &1.4	&-1.3\\
            Random		&50	     &77 &31.28            &3.1	&-2.8\\
			\bottomrule
	\end{tabular}}
	\caption{Base-case Problem: Increased Storage Capacity}
	\label{tbl:Capacity-NoStop3}
\end{table}%

The results of Table~\ref{tbl:Capacity-NoStop3} indicate that a 25\% increase of storage capacity led to 1.1 to 2.5\% increase of reactor's utilization and 1.1 to 2.4\% decrease of costs. A 50\% increase led to 1.1 to 5\% increase of reactor's utilization and 1.1 to 4.5\% decrease of costs. The greatest improvements are observed when every bale processed has high moisture level, and when long sequences are processed. In these problems, the maximum inventory  reached the storage capacity (see Tables \ref{tbl:Cost_base} and \ref{tbl:Cost_sequences}). Thus, by increasing storage, additional inventory accumulated, which led to increased utilization of the reactor. The problems with short sequences show the least improvements since the maximum inventory level was lower than capacity (see Tables \ref{tbl:Cost_base} and \ref{tbl:Cost_sequences}). Thus, increasing storage capacity has  minimal impact on reactor utilization.



\noindent{\bf Short-duration Failures:}
Table~\ref{tbl:Capacity-ShortStop3} summarizes the impact of increased storage capacity on costs and reactor utilization for the cases when reactor's capacity is low and high.
The last two columns of this table present the change in these performance measures as compared to the results of the problem under short duration equipment failures (Table \ref{tbl:Cost_shortfail}).

The results of Table~\ref{tbl:Capacity-ShortStop3} indicate that  a 25\% increase of storage capacity led to 1.2 to 2.7\% increase of reactor’s utilization and 1.2 to 2.6\% decrease of costs.  A 50\% increase led to 2.5 to 4.4\% increase of reactor’s utilization and 2.3 to 4.0\% decrease of costs.  The greatest improvements are observed when  every  bale  processed  has  high  moisture  level,  and  when  long  sequences  are  processed. The increase of storage capacity has a greater impact on reducing costs and increasing reactor's utilization as compared to the base-case problem (Table \ref{tbl:Capacity-NoStop3}). This is because the increase of storage capacity  allows for additional accumulation of inventory which is used during equipment failures.

\begin{table}[htbp]
\centering
	\scriptsize
	{%
		\begin{tabular}{cc|cc|cc}
			\toprule
			\textbf{Biomass Feeding} &\textbf{Inventory Capacity}&\textbf{Reactor}& \textbf{Total}&\textbf{Reactor}& \textbf{Total}\\
			\textbf{Pattern} &\textbf{Increase}&\textbf{Utilization}&\textbf{Cost}&\textbf{Utilization $\Delta$}&\textbf{Cost $\Delta$}\\
		    & $(\%)$ & $(\%)$ & $(\$/dt)$& $(\%)$ & $(\%)$\\
			\hline\\
			\multicolumn{6}{c}{\bf Low Reactor Capacity}\\
		\hline
             High  & 25 & 75 & 58.68& 1.8 & -1.7\\
             High  & 50 & 78 & 57.30 & 4.4 & -4.0\\
            \hline\\
			\multicolumn{6}{c}{\bf High Reactor Capacity}\\
			\hline
			High		&25	 &43&57.89&2.7	&-2.6  \\
             High		&50	 &44 & 57.17   &4.1    &-3.8  \\
              \cmidrule(lr){1-1}
			\cmidrule(lr){2-2}
			\cmidrule(lr){3-3}
			\cmidrule(lr){4-4}
			\cmidrule(lr){5-5}
			\cmidrule(lr){6-6}
            Long		&25	  &71&34.11  &1.9    &-1.7  \\
            Long		&50	  &72& 33.56&3.7    &-3.3  \\
             \cmidrule(lr){1-1}
			\cmidrule(lr){2-2}
			\cmidrule(lr){3-3}
			\cmidrule(lr){4-4}
			\cmidrule(lr){5-5}
			\cmidrule(lr){6-6}
            Short		&25	  &75&32.00  &1.2    &-1.2  \\
            Short		&50	  &76&31.64  &2.5    &-2.3  \\
             \cmidrule(lr){1-1}
			\cmidrule(lr){2-2}
			\cmidrule(lr){3-3}
			\cmidrule(lr){4-4}
			\cmidrule(lr){5-5}
			\cmidrule(lr){6-6}
            Random		&25	    &75& 32.11&2.0	&-1.9  \\
            Random		&50	    &76&31.72 &3.4	&-3.1  \\
			\bottomrule
	\end{tabular}}
	\caption{Short-Duration Failures: Increased Storage Capacity}
	\label{tbl:Capacity-ShortStop3}
\end{table}%

\noindent{\bf Long-duration Failures:} Table~\ref{tbl:Capacity-LongStop3} summarizes the impact of increased storage capacity on costs and reactor utilization for the cases when  reactor's capacity is low and high. The last two columns of this table present the change in these performance measures as compared to the results of the problem under long-duration equipment failures (Table~\ref{tbl:Cost_longfail}).

The results of Table~\ref{tbl:Capacity-LongStop3} indicate that  a 25\% increase of storage capacity led to 3.9 to 33.8\% increase of reactor’s utilization and 3.7 to 24.9\% decrease of costs. A 50\% increase led to 3.9 to 51.7\% increase of reactor’s utilization and 3.7 to 33.8\% decrease of costs. An increase of storage capacity has the greatest impact on reducing costs and increasing reactor utilization when the system experiences long-duration failures. The results of Table~\ref{tbl:Cost_longfail} indicate that the maximum inventory reached storage capacity in all problems solved. Thus, increasing the storage capacity leads to additional accumulation of the inventory which is used to maintain a continuous flow of biomass to the reactor during failures.    
%
\begin{table}[h]
\centering
	\scriptsize
	{%
		\begin{tabular}{cc|cc|cc}
			\toprule
			\textbf{Biomass Feeding} &\textbf{Inventory Capacity}&\textbf{Reactor}& \textbf{Total}&\textbf{Reactor}& \textbf{Total}\\
			\textbf{Pattern} &\textbf{Increase}&\textbf{Utilization}&\textbf{Cost}&\textbf{Utilization $\Delta$}&\textbf{Cost $\Delta$}\\
		    & $(\%)$ & $(\%)$ & $(\$/dt)$& $(\%)$ & $(\%)$\\
			\hline\\
			\multicolumn{6}{c}{\bf Low Reactor Capacity}\\
			\hline
			Low	    &25	  & 82 & 51.03  &3.9	&-3.7  \\
            Low	    &50	  &82 &51.02  &3.9	&-3.7  \\
            \cmidrule(lr){1-1}
			\cmidrule(lr){2-2}
			\cmidrule(lr){3-3}
			\cmidrule(lr){4-4}
			\cmidrule(lr){5-5}
			\cmidrule(lr){6-6}
            High	&25	  & 51 & 87.02  &14.4	&-12.5  \\
            High	&50	  &56 &79.72  &25.0   &-19.9  \\
            \cmidrule(lr){1-1}
			\cmidrule(lr){2-2}
			\cmidrule(lr){3-3}
			\cmidrule(lr){4-4}
			\cmidrule(lr){5-5}
			\cmidrule(lr){6-6}
            Long	&25	 &70&61.26   &14.1	    &-12.3  \\
            Long	&50	  &76&56.64  &23.6	    &-18.9 \\
            \cmidrule(lr){1-1}
			\cmidrule(lr){2-2}
			\cmidrule(lr){3-3}
			\cmidrule(lr){4-4}
			\cmidrule(lr){5-5}
			\cmidrule(lr){6-6}
            Short	&25	 &70 &60.50   &11.4	    &-10.3  \\
            Short	&50	 &77&55.46   &21.6	    &-17.8  \\
            \cmidrule(lr){1-1}
			\cmidrule(lr){2-2}
			\cmidrule(lr){3-3}
			\cmidrule(lr){4-4}
			\cmidrule(lr){5-5}
			\cmidrule(lr){6-6}
            Random	&25	 &65 &65.56   &15.8	    &-13.7  \\
            Random	&50	 &72&59.46   &27.7	    &-21.7  \\
			\hline \\
            \multicolumn{6}{c}{\bf High Reactor Capacity}\\
			\hline
			Low	    &25	 &61 & 39.05   &17.1	    &-14.5  \\
            Low	    &50	 &68 &35.10    &30.3	    &-23.2  \\
            \cmidrule(lr){1-1}
			\cmidrule(lr){2-2}
			\cmidrule(lr){3-3}
			\cmidrule(lr){4-4}
			\cmidrule(lr){5-5}
			\cmidrule(lr){6-6}
            High	&25	 &29 &86.32   &15.2	    &-13.1  \\
            High	&50	 &32 &78.50   &26.8	    &-21.0  \\
            \cmidrule(lr){1-1}
			\cmidrule(lr){2-2}
			\cmidrule(lr){3-3}
			\cmidrule(lr){4-4}
			\cmidrule(lr){5-5}
			\cmidrule(lr){6-6}
            Long	&25	  &45 &54.24  &17.6	    &-14.8  \\
            Long	&50	  &50 &48.04  &32.8	    &-24.6  \\
            \cmidrule(lr){1-1}
			\cmidrule(lr){2-2}
			\cmidrule(lr){3-3}
			\cmidrule(lr){4-4}
			\cmidrule(lr){5-5}
			\cmidrule(lr){6-6}
            Short	&25	 &46 &52.37   &14.2	    &-12.4  \\
            Short	&50	  &51 &47.51  &26.0	    &-20.5  \\
            \cmidrule(lr){1-1}
			\cmidrule(lr){2-2}
			\cmidrule(lr){3-3}
			\cmidrule(lr){4-4}
			\cmidrule(lr){5-5}
			\cmidrule(lr){6-6}
            Random	&25	 &42 &57.64   &33.8	    &-24.9  \\
            Random	&50	  &48 &50.83  &51.7	    &-33.8  \\
            \bottomrule
	\end{tabular}}
	\caption{Long-Duration Failures:  Increased Storage Capacity}
	\label{tbl:Capacity-LongStop3}
\end{table}%

\subsection{Sensitivity Analysis: Biomass  Particle Size}\label{sec:particleSize}
In this section, we analyze the impact of particle size on the performance of the system. Regression equations \eqref{reg1} and \eqref{reg2} show the relationship between particle size distribution (represented by $\rho^{50}$) 
and biomass density. Biomass density impacts the weight of biomass that can be stored in the metering bin. Particle size also impacts the separation process after the first grinder. Our experiments focus on evaluating the impact that ($i$)  particle size ($\rho^{50}$), and ($ii$) particle uniformity ($\frac{\rho^{90}}{\rho^{10}}$) have on the performance of the system. 

 Table~\ref{tab:PSD} in Appendix C summarizes the distribution of particle size for the base-case problem. In the following analysis we use as a reference this problem for the case when low moisture biomass is processed and reactor's capacity is high. Similar observations are made for problems that use other biomass feeding patterns. 
 %

\begin{figure}[h]
    \centering
    \includegraphics[width=0.47\textwidth]{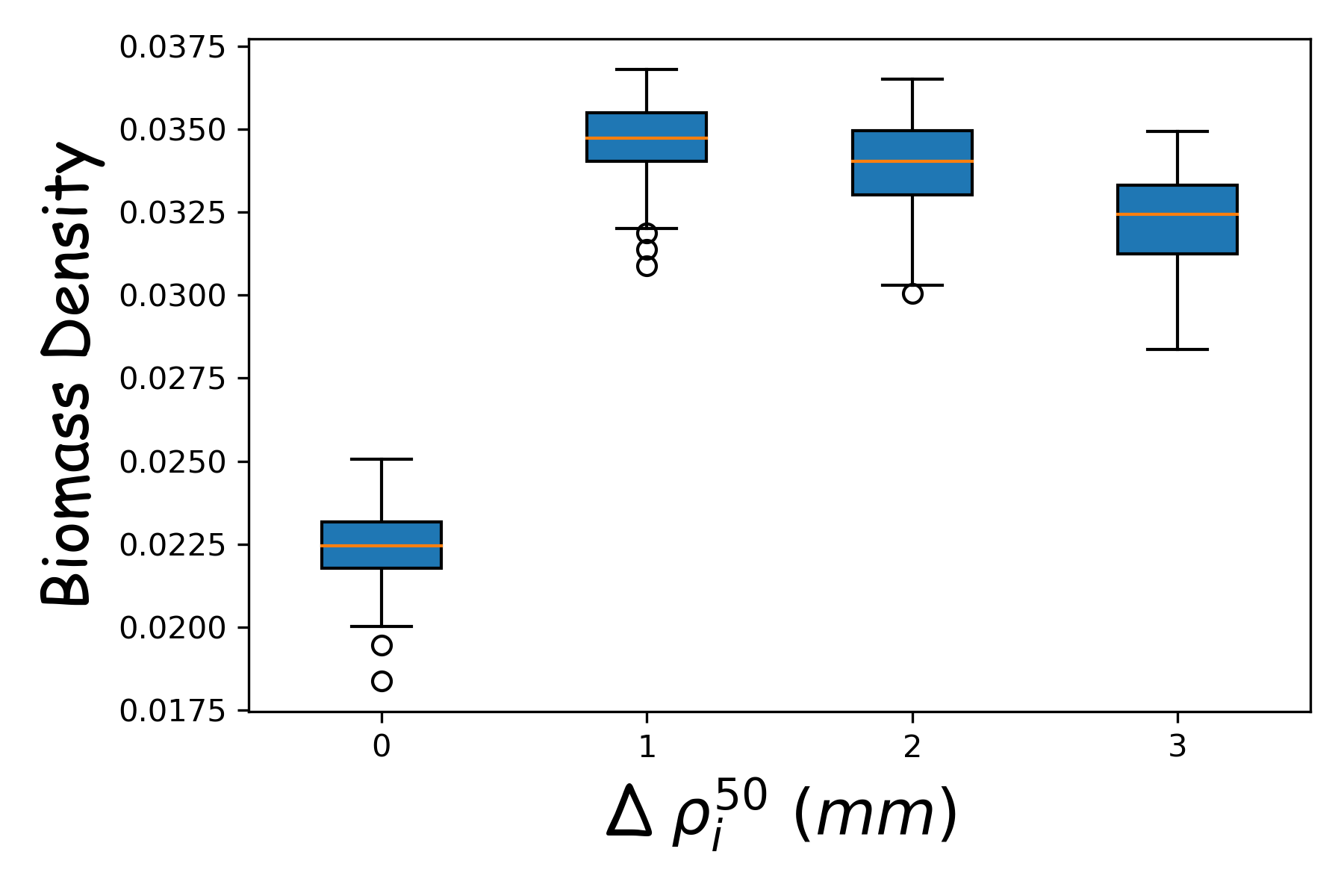}
    \hspace{0.01cm}
    \includegraphics[width=0.47\textwidth]{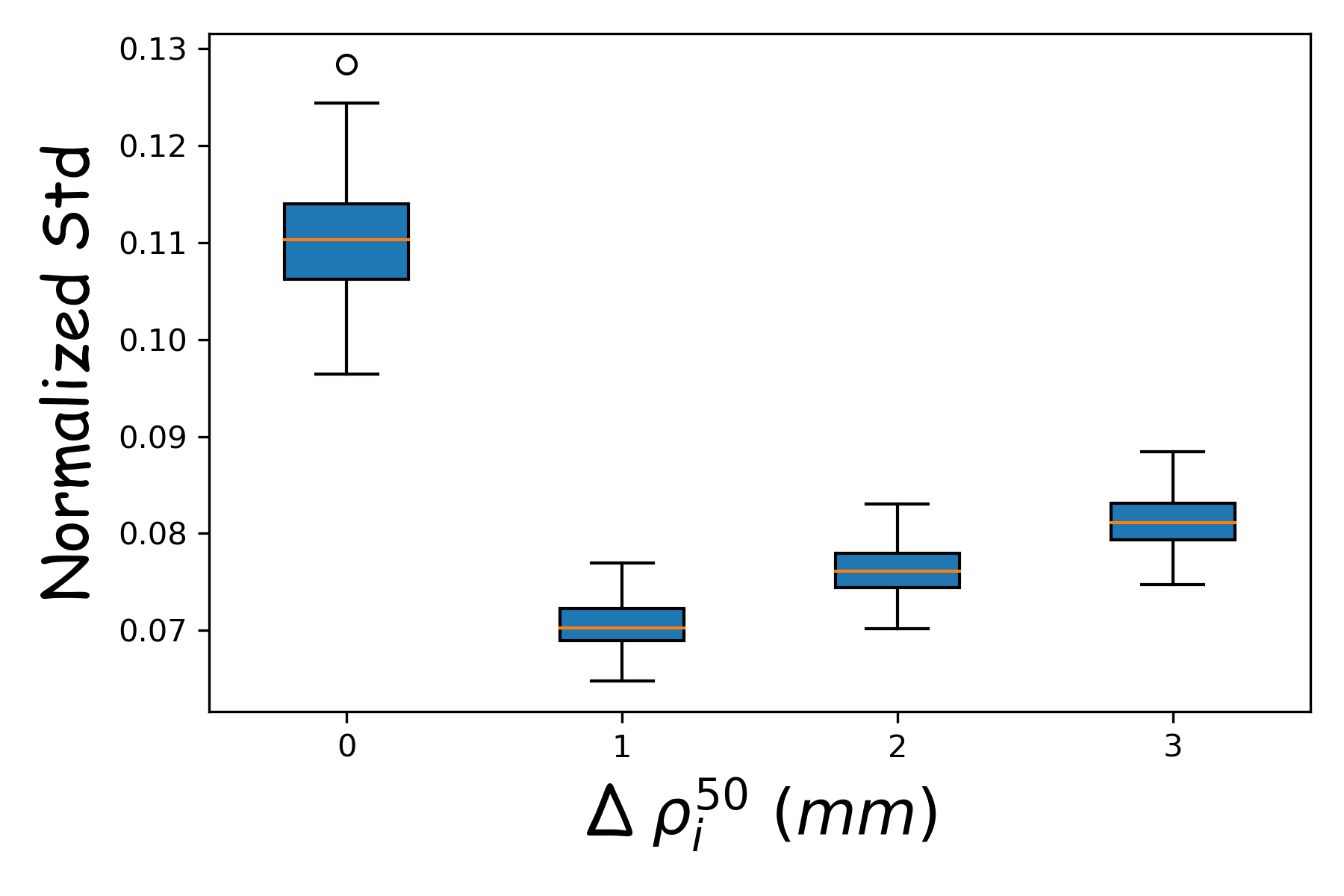}
    \caption{Average and Standard Deviation of Biomass in Metering Bin.}\label{fig:NStdDen_g1D50}
\end{figure}

\noindent {\bf Primary Grinder:} Table~\ref{tbl:Cost_g1D50} summarizes the change of total cost and reactor's utilization due to increases of particle size. The distribution of particle size can be controlled by changing the rotational speed of mills in the primary grinder. Figure \ref{fig:NStdDen_g1D50} summarizes the distribution of the mean particle size in the metering bin, and the corresponding standard deviation. 

One would expect that by increasing particle size, the density of biomass in the metering bin would decrease. Different from this intuition, an initial increase of particle size leads to an increase of biomass density in the metering bin. This is mainly because, a larger proportion of biomass will need to be reprocessed in the secondary grinder to reduce particle size. Reprocessing leads to particles of smaller size (as compared to processing only on the primary grinder) in the metering bin, which leads to increased biomass density and decreased volume. As a result, additional biomass can be stored in the metering bin, which enables the system to maintain a continuous flow of biomass to the reactor. This leads to an increase of reactor's utilization and a decrease of costs. 

However, further increases of particle size decreases biomass density, what leads to an decrease of the weight of biomass flow in the equipment that feed the secondary grinder. Specifically, the weight of biomass that is moved via drag chain (DC 5 \& DC6) and screw (SC6) conveyors decreases. For example, when mean particle size increases by 3mm, SC6 can transport 1.0 to 6.2 $ton/ hr$. This limits the flow of biomass to the reactor.         

\begin{table}[htbp]
\centering
\scriptsize
	{%
\centering
		\begin{tabular}{cccc|cc|cc|cc}
			\toprule
			$\Delta \rho^{50}$&\textbf{Bypass }&\textbf{$\Delta$ Reactor} &\textbf{$\Delta$ Total} & \multicolumn{2}{c|}{\textbf{D5 Capacity}}& \multicolumn{2}{c|}{\textbf{D6 Capacity}}& \multicolumn{2}{c}{\textbf{SC6 Capacity}}\\
			&\textbf{Ratio}&\textbf{Utilization}&\textbf{Cost} & \bf Min & \bf Max& \bf Min & \bf Max& \bf Min & \bf Max\\
		    $(mm)$& $(\%)$ & $(\%)$ & $(\%)$ & $(ton/hr)$& $(ton/hr)$& $(ton/hr)$& $(ton/hr)$& $(ton/hr)$& $(ton/hr)$\\
			\cmidrule(lr){1-1}
			\cmidrule(lr){2-2}
			\cmidrule(lr){3-3}
			\cmidrule(lr){4-4}\cmidrule(lr){5-5}\cmidrule(lr){6-6}\cmidrule(lr){7-7}\cmidrule(lr){8-8}\cmidrule(lr){9-9}\cmidrule(lr){10-10}
            0 & 87 & 0.0 & 0.0 & 11.5 & 21.2 &11.5&	21.2&6.2	&11.3\\
            +1&	58 & +6.0 & 	 -5.6 & 7.7	&17.3 &7.7&	17.3&4.1	&9.3\\
            +2&	53& +4.7&    -4.5  &5.8	&15.4 &5.8	&15.4&3.1 & 8.2\\ 
            +3&	51 & -3.1&	 +3.3 &  1.9&	11.5 &1.9&	11.5 &1.0 &	6.2\\
			\bottomrule
	\end{tabular}}
	\caption{Particle Size: Primary Grinder}
	\label{tbl:Cost_g1D50}
\end{table}
\begin{table}[htbp]
\centering
\scriptsize
	{%
	\begin{minipage}{0.45\textwidth}
	  \centering
    \begin{tabular}{cccc}
		\toprule
		\textbf{$\Delta \frac{\rho^{90}}{\rho^{10}}$}&\textbf{Bypass} & \textbf{$\Delta$ Reactor} &\textbf{$\Delta$ Total}\\
		&\textbf{Ratio}&\textbf{Utilization}&\textbf{Cost} \\
	    & $(\%)$ & $(\%)$ & $(\%)$\\
		\cmidrule(lr){1-1}
		\cmidrule(lr){2-2}
		\cmidrule(lr){3-3}
		\cmidrule(lr){4-4}
        -1&	91 & -2.2&	 2.3   \\
        -2& 95&	-3.4&	 3.5   \\ 
        -3& 100&	-1.8&	 1.8   \\
		\bottomrule
	\end{tabular}
	\caption{Particle Uniformity: Primary Grinder}
	\label{tbl:Cost_g1D90D10}
	\end{minipage}
	\begin{minipage}{0.45\textwidth}
	\centering
	\begin{tabular}{ccc}
			\toprule
			$\Delta \rho^{50}$&\textbf{$\Delta$ Reactor} &\textbf{$\Delta$ Total}\\
			& \textbf{Utilization}&\textbf{Cost} \\
		    $(mm)$& $(\%)$ & $(\%)$\\
			\cmidrule(lr){1-1}
			\cmidrule(lr){2-2}
			\cmidrule(lr){3-3}
           + 0.1&   + 0.8&	 -0.8   \\
           + 0.2&	+1.5&	 -1.4   \\ 
           + 0.3&  +2.3&	 -2.3   \\ 
			\bottomrule
	\end{tabular}
	\caption{Particle Size: Secondary Grinder}
	\label{tbl:Cost_g2D50}
	\end{minipage}
	
	}
\end{table}%

\begin{figure}
    \centering
     \includegraphics[scale=0.45]{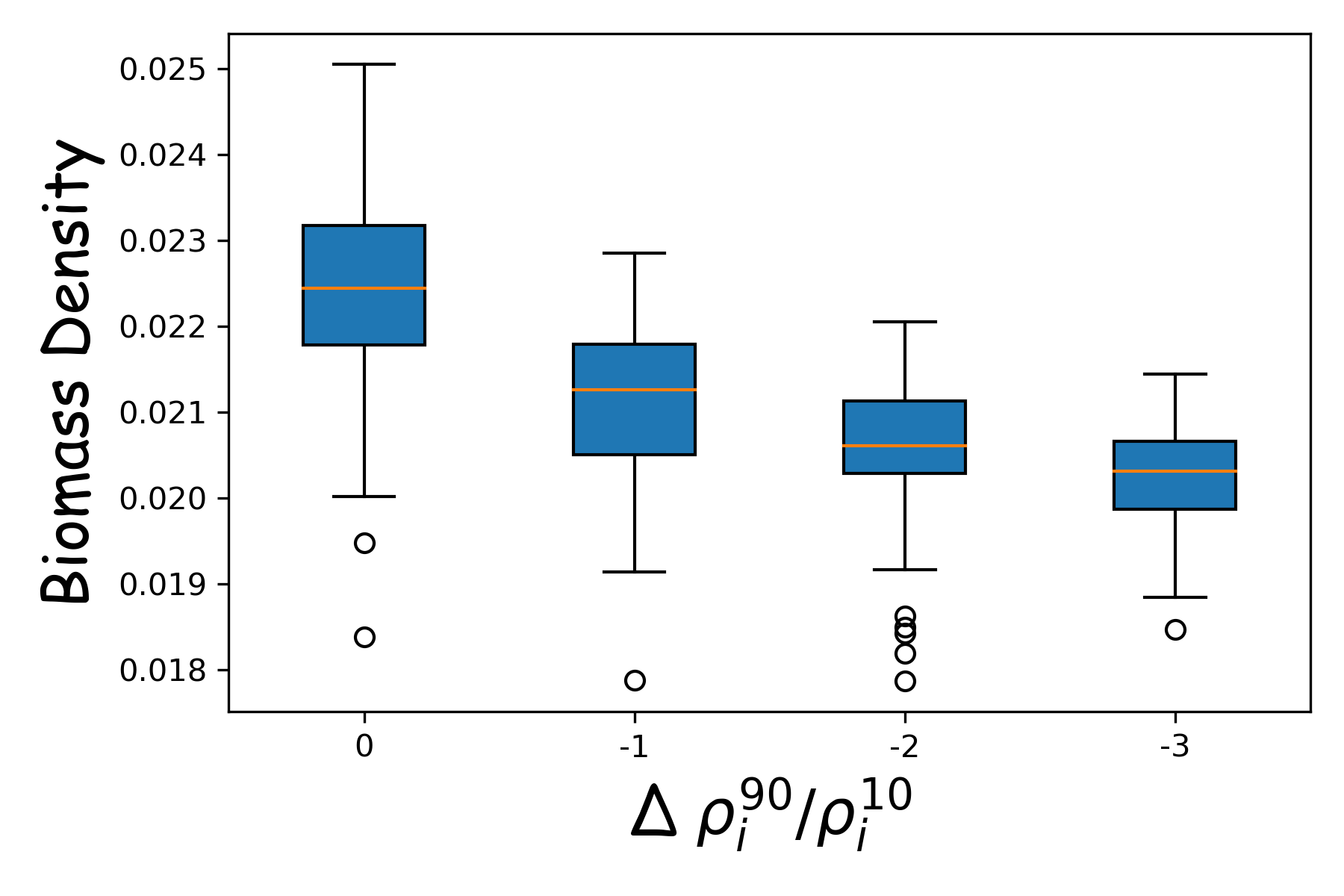}
     \caption{Average Biomass Density}
     \label{fig:AvgDen_g1D9010}
    \includegraphics[width=0.45\textwidth]{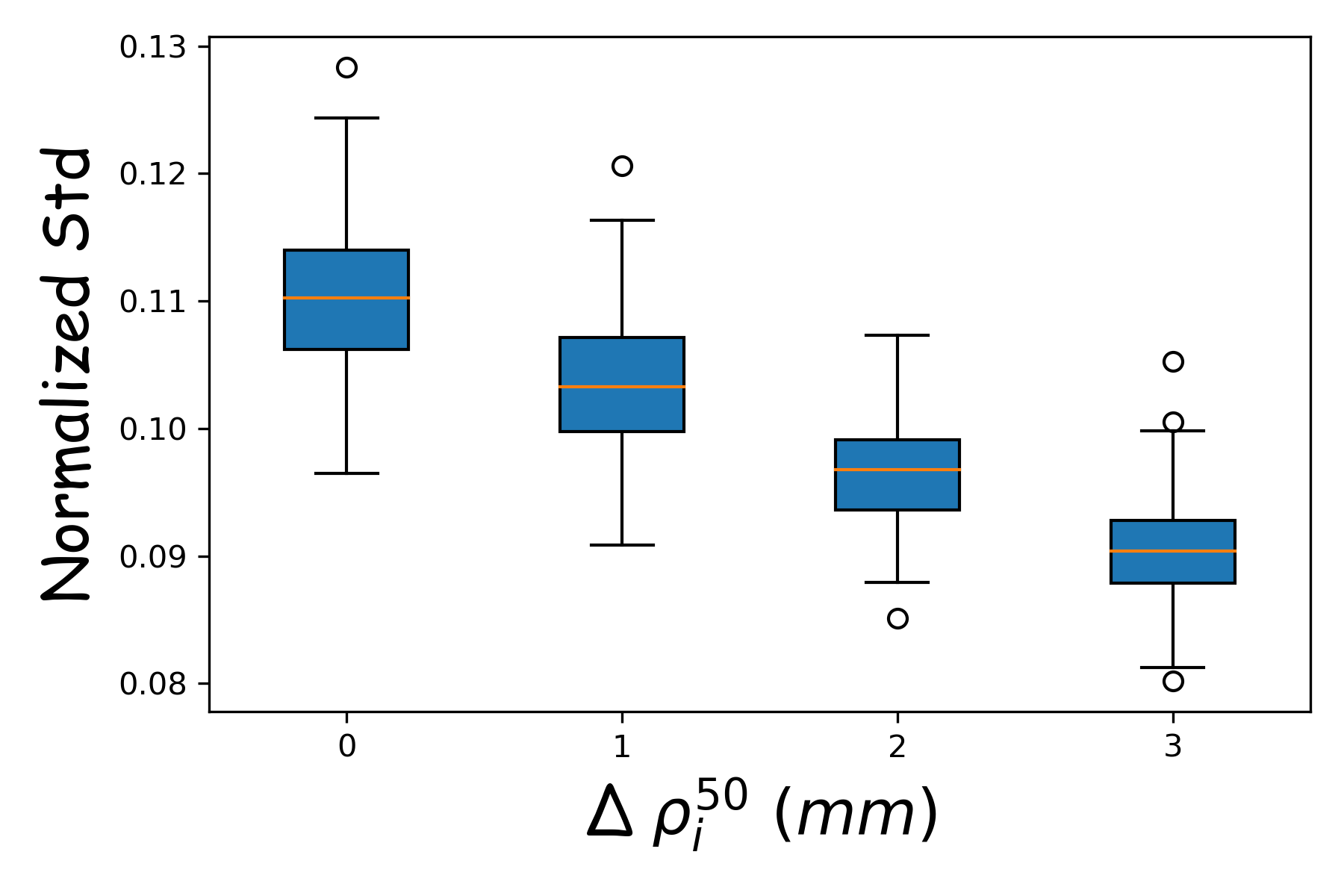}
    \caption{Normalized Standard Deviation of Biomass Density in the Metering Bin}
    \label{fig:Nor-Std_g2D50}
\end{figure}

Table~\ref{tbl:Cost_g1D90D10} summarizes the results of changing the particle size uniformity on reactor utilization and costs. Notice that, by decreasing $\frac{\rho^{90}}{\rho^{10}}$, we  increase particle size uniformity, which means, the distribution of particle size gets close to the mean value. This leads to an increased bypass ratio, which is calculated via equation \eqref{eq:bypassRatio}. Figure~\ref{fig:AvgDen_g1D9010} indicates a decrease of biomass density in the metering bin. These changes lead to decreases of reactor utilization and increases cost. 

\noindent {\bf Secondary Grinder:} 
Tables~\ref{tbl:Cost_g2D50} summarizes the impact of changes in particle size after the secondary grinder to the total cost and reactor utilization. Changes in particle size and particle uniformity after the secondary grinder do not impact  the bypass ratio and the flow of equipment upstream. Thus we did not consider changes in particle uniformity after the secondary grinder. Figure \ref{fig:Nor-Std_g2D50} presents the decrease of normalized standard deviation of biomass density in the metering bin. Biomass density in the metering bin follows a similar trend.  The results show that increase of particle size leads to increased utilization of the reactor and decreased total costs. 
\exclude{
\begin{table}[htbp]
\centering

\centering \scriptsize
	{%
\centering

		\begin{tabular}{ccc}
			\toprule
			$\Delta \rho^{50}$&\textbf{$\Delta$ Reactor} &\textbf{$\Delta$ Total}\\
			& \textbf{Utilization}&\textbf{Cost} \\
		    $(mm)$& $(\%)$ & $(\%)$\\
			\cmidrule(lr){1-1}
			\cmidrule(lr){2-2}
			\cmidrule(lr){3-3}
           + 0.1&   + 0.8&	 -0.8   \\
           + 0.2&	+1.5&	 -1.4   \\ 
           + 0.3&  +2.3&	 -2.3   \\ 
			\bottomrule
	\end{tabular}}
	\caption{Particle Size: Secondary Grinder}
	\label{tbl:Cost_g2D50}
\end{table}%
}
Since the particle size of biomass about to be processed in the secondary grinder is smaller than that of biomass processed in the primary grinder, we consider smaller changes of particle size. Thus, the resulting changes to reactor utilization and costs are smaller.

\section{Summary and Conclusion}\label{sec:Conclusion}
The  main  objective  of  this  work  is  to develop analytical models which improve the  reliability  of  biomass  feeding  system  in  a biorefinery. We propose a stochastic programming model that identifies operating conditions of equipment and inventory level to maintain a continuous flow of biomass to the  reactor under uncertainty.  The uncertain problem parameters are moisture content, particle size distribution, and equipment failures.  We use a chance constraint to ensure that a certain reliability level of the reactor is maintained. We propose a Sample Average Approximation of the chance-constraint problem and develop a bisection search-based heuristic to solve this approximation.

This research makes a number of contributions to the literature. {First}, the proposed model is developed using real life data about biomass preprocessing. The model integrates the outcomes of a Discrete Element Model to an optimization model. As a result, we expect that practitioners will find the results of this study useful and applicable. 

{Second}, we make a number of important observations. We observe that ($i$) \emph{biomass characteristics} impact reactor utilization and costs. The system performs best when processing biomass of low moisture and performs worst when processing biomass of high moisture content. Changes of particle size distribution also impact the performance of the system by increasing reactor's utilization by as much as 6\% and reducing costs by as much as 5.6\%. We observe that ($ii$) \emph{sequencing bales} based on moisture level improves the use of resources, leading to up to 2.6\% increase of reactor's utilization and 1.96\% reduction of costs. We observed that short sequences of low, medium and high moisture bales perform best. Processing low and medium moisture bales builds up the inventory required to maintain a continuous flow of biomass when processing high moisture bales since their processing times are long. We observe that ($iii$) short equipment failures reduce reactor's utilization by up to 4.0\%, increase costs by up to 4.2\%, and increase inventory level by up to 55.8\%. Long duration failures reduce reactor's utilization by up to 58.1\%, increase costs by up to 139.2\%, increase average inventory by up to 225.8\%. The limitations of storage capacity negatively impacted the performance of the system during long duration failures. Thus,  increasing storage capacity and processing short sequences of low, medium and high moisture biomass are strategies that can be use to mitigate the negative impacts of equipment failure. 

We identify several research directions that are worth exploring in the future. First, the model can be extended to consider adaptive changes to the operating conditions of equipment in every period based on the specific characteristics of biomass. This will require the development of a multi-stage stochastic programming model, which utilizes sensor-based biomass data collected during the process to make real-time adaptive decisions. Second, the proposed model can be integrated with a macro-level supply chain model which identifies supplier of a biorefinery based on the moisture level of biomass, selling price, transportation cost, etc. The model will capture the trade-offs between biomass selling price and the cost of processing biomass. Finally, the model can be extended to identify what sequence of biomass bales  has the greatest impact on the performance of the system. The proposed model only sequences bales based on moisture level, however, additional biomass characteristics can be used. In this case, the complexity of the problem increases, which would motivate the development of efficient solution algorithms.

\bibliographystyle{plain}
\bibliography{msp-forecast-error}

\begin{thebibliography}{10}

\bibitem{abdelaziz2012solution}
F.~B. Abdelaziz.
\newblock Solution approaches for the multiobjective stochastic programming.
\newblock {\em European Journal of Operational Research}, 216(1):1--16, 2012.

\bibitem{Alsina2018}
E.F. Alsina, M.~Chica, and K.~Trawinski.
\newblock On the use of machine learning methods to predict component
  reliability from data-driven industrial case studies.
\newblock {\em The International Journal of Advanced Manufacturing Technology},
  94:2419 – 2433, 2018.

\bibitem{atlason2008optimizing}
J.~Atlason, M.A. Epelman, and S.G. Henderson.
\newblock Optimizing call center staffing using simulation and analytic center
  cutting-plane methods.
\newblock {\em Management Science}, 54(2):295--309, 2008.

\bibitem{Basciftci2018}
B.~{Basciftci}, S.~{Ahmed}, N.~Z. {Gebraeel}, and M.~{Yildirim}.
\newblock Stochastic optimization of maintenance and operations schedules under
  unexpected failures.
\newblock {\em IEEE Transactions on Power Systems}, 33(6):6755--6765, 2018.

\bibitem{Bellman34}
R.~Bellman.
\newblock Dynamic programming.
\newblock {\em Science}, 153(3731):34--37, 1966.

\bibitem{Birge_Louveaux_1997}
J.R. Birge and F.~Louveaux.
\newblock {\em Introduction to Stochastic Programming}.
\newblock Springer, New York, 1997.

\bibitem{Cafferty2013BLM}
K.~Cafferty, D.~Muth, J.~Jacobson, and K.~Bryden.
\newblock Model based biomass system design of feedstock supply systems for
  bioenergy production.
\newblock volume~2, 08 2013.

\bibitem{charnes1955optimal}
Abraham Charnes, William~W Cooper, and Robert~O Ferguson.
\newblock Optimal estimation of executive compensation by linear programming.
\newblock {\em Management science}, 1(2):138--151, 1955.

\bibitem{Chen_Onal_2014}
X.~Chen and H.~Önal.
\newblock An economic analysis of the future u.s. biofuel industry, facility
  location, and supply chain network.
\newblock {\em Transportation Science}, 48(4):575--591, 2014.

\bibitem{DavidCoit_reliability_opt2018}
D.W. Coit and E.~Zio.
\newblock The evolution of system reliability optimization.
\newblock {\em Reliability Engineering and System Safety}, 2018.

\bibitem{Crawford_2016_Flowability_CornStover}
N.~Crawford, N.~Nagle, D.~Sievers, and J.~Stickel.
\newblock The effects of physical and chemical preprocessing on the flowability
  of corn stover.
\newblock {\em Biomass and Bioenergy}, 85:126--134, 02 2016.

\bibitem{DEM1979}
P.~A. Cundall and O.~D.~L. Strack.
\newblock A discrete numerical model for granular assemblies.
\newblock {\em Géotechnique}, 29(1):47--65, 1979.

\bibitem{DAI2012716}
J.~Dai, H.~Cui, and J.R. Grace.
\newblock Biomass feeding for thermochemical reactors.
\newblock {\em Progress in Energy and Combustion Science}, 38(5):716 -- 736,
  2012.

\bibitem{DOMAC2005Bioenergy}
J.~Domac, K.~Richards, and S.~Risovic.
\newblock Socio-economic drivers in implementing bioenergy projects.
\newblock {\em Biomass and Bioenergy}, 28(2):97 -- 106, 2005.

\bibitem{EKSIOGLU20091342}
S.D. Ekşioğlu, A.~Acharya, L.E. Leightley, and S.~Arora.
\newblock Analyzing the design and management of biomass-to-biorefinery supply
  chain.
\newblock {\em Computers \& Industrial Engineering}, 57(4):1342 -- 1352, 2009.

\bibitem{Fyffe1986Reliability}
D.~E. {Fyffe}, W.~W. {Hines}, and N.~K. {Lee}.
\newblock System reliability allocation and a computational algorithm.
\newblock {\em IEEE Transactions on Reliability}, R-17(2):64--69, 1968.

\bibitem{Ghare1969Reliability}
P.~M. Ghare and R.~E. Taylor.
\newblock Optimal redundancy for reliability in series systems.
\newblock {\em Operations Research}, 17(5):838--847, 1969.

\bibitem{guo2020discrete}
Y.~Guo, Q.~Chen, Y.~Xia, T.~Westover, S.D. Eksioglu, and M.~Roni.
\newblock Discrete element modeling of switchgrass particles under compression
  and rotational shear.
\newblock {\em Biomass and Bioenergy}, 141:105649, 2020.

\bibitem{HANSEN2019BiomassSC}
J.K. Hansen, M.S. Roni, S.K. Nair, D.S. Hartley, L.M. Griffel, V.~Vazhnik, and
  S.~Mamun.
\newblock Setting a baseline for integrated landscape design: Cost and risk
  assessment in herbaceous feedstock supply chains.
\newblock {\em Biomass and Bioenergy}, 130:105388, 2019.

\bibitem{Hohner2012}
D.~Höhner, S.~Wirtz, and V.~Scherer.
\newblock A numerical study on the influence of particle shape on hopper
  discharge within the polyhedral and multi-sphere discrete element method.
\newblock {\em Powder Technology}, 226:16–28, 08 2012.

\bibitem{osti_1369631}
J.J. Jacobson, P.~Lamers, M.S. Roni, K.G. Cafferty, K.L. Kenney, B.M. Heath,
  and J.K. Hansen.
\newblock Techno-economic analysis of a biomass depot.
\newblock 10 2014.

\bibitem{osti_1133890}
K.L. Kenney, K.G. Cafferty, J.J. Jacobson, I.J. Bonner, G.L. Gresham, J.R.
  Hess, W.A. Smith, D.N. Thompson, V.S. Thompson, J.S. Tumuluru, and N.~Yancey.
\newblock Feedstock supply system design and economics for conversion of
  lignocellulosic biomass to hydrocarbon fuels conversion pathway: Fast
  pyrolysis and hydrotreating bio-oil pathway "the 2017 design case".
\newblock 1 2014.

\bibitem{BeyondBeta}
S.~Kotz and J.R. van Dorp.
\newblock {\em Beyond Beta}.
\newblock WORLD SCIENTIFIC, 2004.

\bibitem{Li_ChanceConstr_Reliability_2008}
P.~Li, H.~Arellano-Garcia, and G.~Wozny.
\newblock Chance constrained programming approach to process optimization under
  uncertainty.
\newblock {\em Computers and Chemical Engineering}, 32(1):25 -- 45, 2008.
\newblock Process Systems Engineering: Contributions on the State-of-the-Art.

\bibitem{Luedtke_Ahmed_2008}
J.~Luedtke and S.~Ahmed.
\newblock A sample approximation approach for optimization with probabilistic
  constraints.
\newblock {\em {SIAM} Journal on Optimization}, 19:674--699, 2008.

\bibitem{Memisoglu2015}
G.~Memişoğlu and H.~Üster.
\newblock Integrated bioenergy supply chain network planning problem.
\newblock {\em Transportation Science}, 50(1):35--56, 2016.

\bibitem{MENG2020773}
Z.~Meng, Z.~Zhang, and H.~Zhou.
\newblock A novel experimental data-driven exponential convex model for
  reliability assessment with uncertain-but-bounded parameters.
\newblock {\em Applied Mathematical Modelling}, 77:773 -- 787, 2020.

\bibitem{Misra_IP-Reliability_1991}
K.B. Misra.
\newblock An algorithm to solve integer programming problems: An efficient tool
  for reliability design.
\newblock {\em Microelectronics Reliability}, 31(2):285 -- 294, 1991.

\bibitem{NUMBI20161653}
B.P. Numbi and X.~Xia.
\newblock Optimal energy control of a crushing process based on vertical shaft
  impactor.
\newblock {\em Applied Energy}, 162:1653 -- 1661, 2016.

\bibitem{OREFICE2017347}
L.~Orefice and J.G. Khinast.
\newblock Dem study of granular transport in partially filled horizontal screw
  conveyors.
\newblock {\em Powder Technology}, 305:347 -- 356, 2017.

\bibitem{Pagnoncelli_Ahmed_Shapiro_2009}
B.~Pagnoncelli, S.~Ahmed, and A.~Shapiro.
\newblock Sample average approximation method for chance constrained
  programming: Theory and applications.
\newblock {\em J. Optim. Theory Appl}, 142:399--416, 2009.

\bibitem{Painton1995PCReliability}
L.~{Painton} and J.~{Campbell}.
\newblock Genetic algorithms in optimization of system reliability.
\newblock {\em IEEE Transactions on Reliability}, 44(2):172--178, 1995.

\bibitem{Viet2012BiorefOpt}
V.~Pham and M.~El-Halwagi.
\newblock Process synthesis and optimization of biorefinery configurations.
\newblock {\em AIChE Journal}, 58(4):1212--1221, 2012.

\bibitem{Prasad_ReliabilityOpt_2000}
V.~R. {Prasad} and W.~{Kuo}.
\newblock Reliability optimization of coherent systems.
\newblock {\em IEEE Transactions on Reliability}, 49(3):323--330, Sep. 2000.

\bibitem{RONI2014115}
M.S. Roni, S.D. Eksioglu, E.~Searcy, and K.~Jha.
\newblock A supply chain network design model for biomass co-firing in
  coal-fired power plants.
\newblock {\em Transportation Research Part E: Logistics and Transportation
  Review}, 61:115 -- 134, 2014.

\bibitem{SCHERER2016896}
V.~Scherer, M.~Mönnigmann, M.O. Berner, and F.~Sudbrock.
\newblock Coupled dem–cfd simulation of drying wood chips in a rotary drum
  – baffle design and model reduction.
\newblock {\em Fuel}, 184:896 -- 904, 2016.

\bibitem{Sims2003Bioenergy}
R.~Sims.
\newblock Bioenergy to mitigate for climate change and meet the needs of
  society, the economy and the environment.
\newblock {\em Mitigation and Adaptation Strategies for Global Change},
  8(4):349--370, 2003.

\bibitem{Thomopoulos2017}
N.T. Thomopoulos.
\newblock {\em Triangular}, pages 107--112.
\newblock Springer International Publishing, Cham, 2017.

\bibitem{DOE2016}
DOE U.S. Department~of Energy.
\newblock Biorefinery optimization workshop summary report, 2016.

\bibitem{XIA20191}
Y.~Xia, Z.~Lai, T.~Westover, J.~Klinger, H.~Huang, and Q.~Chen.
\newblock Discrete element modeling of deformable pinewood chips in cyclic
  loading test.
\newblock {\em Powder Technology}, 345:1 -- 14, 2019.

\bibitem{yancey2015size}
N.~Yancey and T.~JayaShankar.
\newblock Size reduction, drying and densification of high moisture biomass.
\newblock {\em Quarterly Progress Report. Idaho Fall, Idaho, USA: Idaho
  National Laboratory}, 2015.

\bibitem{YANG2010671}
B.~Yang, X.~Li, M.~Xie, and F.~Tan.
\newblock A generic data-driven software reliability model with model mining
  technique.
\newblock {\em Reliability Engineering and System Safety}, 95(6):671 -- 678,
  2010.

\bibitem{YANG1999Reliability_GA}
J.-E. Yang, M.-J. Hwang, T.-Y Sung, and Y.~Jin.
\newblock Application of genetic algorithm for reliability allocation in
  nuclear power plants.
\newblock {\em Reliability Engineering $\&$ System Safety}, 65(3):229 -- 238,
  1999.

\bibitem{You2012Biofuel}
F.~You, L.~Tao, D.J. Graziano, and S.W. Snyder.
\newblock Optimal design of sustainable cellulosic biofuel supply chains:
  Multiobjective optimization coupled with life cycle assessment and
  input–output analysis.
\newblock {\em AIChE Journal}, 58(4):1157--1180, 2012.

\bibitem{ZHANG20101929}
S.~Zhang and X.~Xia.
\newblock Optimal control of operation efficiency of belt conveyor systems.
\newblock {\em Applied Energy}, 87(6):1929 -- 1937, 2010.

\bibitem{ZhouCornStover2008}
B.~Zhou, K.~Ileleji, and G.~Ejeta.
\newblock Physical property relationships of bulk corn stover particles.
\newblock {\em Transactions of the ASABE}, 51:581--590, 03 2008.

\bibitem{ZONDERVAN2011BiorefOpt}
E.~Zondervan, M.~Nawaz, A.B. {de Haan}, J.M. Woodley, and R.~Gani.
\newblock Optimal design of a multi-product biorefinery system.
\newblock {\em Computers $\&$ Chemical Engineering}, 35(9):1752 -- 1766, 2011.
\newblock Energy Systems Engineering.

\end{thebibliography}

\pagebreak
\appendix
\section{Detailed Notations}
\begin{table}[htbp]
	\scriptsize
	{%
		\begin{tabular}{ll}
			\toprule
			\multicolumn{2}{l}{\textbf{SETS:}}\\
			 ${\bf M}$ & The set of moisture levels of biomass, ${\bf M} := \{Low, Medium, High\}$. \\
		    $\mathcal{T}$ & The set of time periods in the planning horizon, $\mathcal{T} := \{1,2,\ldots, T\}.$\\
		    ${\bf E}^g$ & The set of grinders.\\
		    ${\bf E}^{mill}$ & The set of pelleting mills.\\
		    ${\bf E}^p$ & The set of processing equipment, ${\bf E}^p := {\bf E}^g \cup {\bf E}^{mill}$. \\
		    ${\bf E}^r$ & The set of transportation equipment.\\
		    ${\bf E}^m$ & The set of storage equipment.\\
		    ${\bf N}$ & The full set of equipment, ${\bf N} := {\bf E}^p \cup {\bf E}^r \cup {\bf E}^m$.\\
		    $\boldsymbol{\delta}^-_i$ & The set of equipment that are connected to the equipment $i \in {\bf N}$. \\
			$\boldsymbol{\mathcal{S}}$ & The set of scenarios, $\boldsymbol{\mathcal{S}}$ $:=\ \{1,2,\dots, S\}$.\\
			\midrule
			\multicolumn{2}{l}{\textbf{PARAMETERS:}}\\
			$T$ & Length of the planning horizon.\\
			$\mathbf{h}$ & Inventory holding costs for storage equipment (in \$/ton).\\
	        $f(\cdot)$ & Energy consumption and operational cost function.\\
			$\kappa_t$ & Moisture level of the biomass bales in time period t ($\kappa_t \in M$).\\	
			$\bar{v}_{i}(\kappa_t)$ & Upper bound of processing speed of equipment $i \in {\bf N}$ in time period $t$.\\ 
			$\bar{\iota}_i(\kappa_0)$ & Initial inventory holding capacity in the storage equipment $i \in {\bf E}^m$.\\
			$\bar{r}$ & Feeding capacity of the reactor (in $dt/t)$.\\
			$S$ & The number of scenarios.\\
			$\epsilon$ & Risk tolerance parameter.\\
			$b$ & Index representing the equipment that feeds the separation unit.\\
			$b_1, b_2$ & Indices representing the first equipment on the secondary grinding branch and bypass branch, respectively.\\
			$k$ & Index representing the last equipment that feeds the reactor.\\
			$\gamma_{0}$ & Cross section area of a bale (in $inch^2$).\\
			$\gamma_i$ & Cross section area for equipment $i \in {\bf E}^r$ and cross section area of the discharge opening for equipment $i \in {\bf E}^p$  (both in $inch^2)$.\\
			$\phi_i$ & Dry matter loss in grinder $i \in {\bf E}^g$.\\
			$\varphi_{it}(\kappa_t)$ & Moisture loss in processing equipment $i \in {\bf E}^p$.\\
			$\hat{t}_i$& Processing time inside the pelleting mill  $i \in \boldsymbol{E^m}$.\\
			$c_i$ & Fixed operational costs of equipment $i \in {\bf N}$ (in $\$$/$hour$).\\
			$e_i$ & Energy consumption of equipment $i \in {\bf N}$  (in $\$$/$hour$).\\
			$\bar{\iota}_i$ & Volumetric capacity of  equipment $i \in {\bf E}^r \cup \boldsymbol{E^m}$ (in $inch^3)$.\\
			$\underline{\iota}_i$ & Volume of the minimum inventory required in the storage equipment $i \in {\bf E}^r$ (in $inch^3)$.\\
			$\bar{x}_i(\kappa)$ & Throughput capacity of equipment $i$ (in $tons/time\ period)$.\\
			\midrule
			\multicolumn{2}{l}{\textbf{RANDOM VARIABLES (SCENARIO-DEPENDENT PARAMETERS):}}\\
			$d_{0ts}$ & Density of a biomass bale in time period $t$ in scenario $s$ (in $tons/inch^3$).\\
			$m_{its}$ & Moisture content of the biomass that flows from the equipment $i \in {\bf E}^g$ in period $t$ in scenario $s$.\\
			$d_{its}$ & Density of biomass flowing from the equipment $i \in {\bf E}^g \cup {\bf E}^m$ in time period $t$ in scenario $s$ (in $tons/inch^3$).\\
			$\rho^j_{its}$ & The $j$-th percentile of the particle size distribution of biomass processed in equipment $i \in {\bf E}^g$ in time period $t$ in scenario $s$.\\
			$\theta_{ts}$ & Proportion of the biomass that bypass the additional grinding in time period $t$ in scenario $s$.\\
			\midrule
			\multicolumn{2}{l}{\textbf{DECISION VARIABLES:}}\\
			$V_{it}$ & Processing speed of the equipment $i \in {\bf N}$ in time period $t$ (in meters or rotations per time period).\\
			$I_{i0}$ & Initial inventory level in the storage equipment $i \in {\bf E}^m$ (in tons). \\
			$I_{its}$ & Inventory level in the storage equipment $i \in {\bf E}^m$ in time period $t$ in scenario $s$ (in tons).\\
			$X_{0ts}$ & Biomass flow to the system (i.e., feeding the primary grinder) in time period $t$ in scenario $s$ (in tons).\\
			$X_{its}$ & Biomass flow from equipment $i \in {\bf N}$ in time period $t$ in scenario $s$ (in tons).\\
			\bottomrule
			\end{tabular}}
\end{table}%
\vfill
\pagebreak
\section{Detailed Formulation}

{\bf Introducing a surrogate linear objective function.}
The true objective of this problem is to minimize the total operational cost per dry tons of biomass processed during the planning horizon, which is computed using unit ``dollar per dry ton'' ($\$/dt$). In this case the energy consumption and operation cost function, $\hat{f}(\cdot)$ is as follows: \vspace{-0.8in}
\begin{equation*}
     \hat{f}(X_{ts}, I_{ts}, \omega_s) = \frac{ \displaystyle \sum_{i \in {\bf N}} T(e_i + c_i)}{\displaystyle \sum_{t\in \mathcal{T}} (1-m_{kts}) X_{kts} }.
\end{equation*}
\vspace{-0.8in}
As a result the true objective that we evaluate is given by: 
\[
\min \ \sum_{i \in {\bf E}^m}h_i I_{i0} 
+   \frac{1}{S}  \frac{\displaystyle \sum_{i\in {\bf N}} T(e_i + c_i) }{\displaystyle \left \lbrack \sum_{s=1}^{S} \sum_{t\in \mathcal{T}} (1-m_{kts}) X_{kts} \right \rbrack }.
\]

However, this results in a non-convex objective function, which makes the problem hard to solve. Notice that, in order to minimize the total cost per dry ton, we need to maximize the amount of biomass processed, while using minimum initial inventory. By this observation, we consider the following surrogate linear objective function, which makes the problem much easier to solve.
\[
\min \ \sum_{i \in E^m}  h_i I_{i0}  - \left( \frac{1}{S} \left \lbrack \sum_{s=1}^{S} \sum_{t\in \mathcal{T}} (1-m_{kts})X_{kts} \right \rbrack  \right).
\]

Note that two terms (i.e., holding cost and average amount of biomass processed) of this surrogate objective function have different units. In order to resolve this we used a weighted sum of the two terms, where the weight of the initial holding cost is much smaller.

Use of this surrogate objective function leads to the following chance-constrained stochastic linear program:
\begin{footnotesize}
\begin{alignat}{2}
(\hat{P}) := \min \ & \sum_{i \in E^m}  h_i I_{i0}  - \left( \frac{1}{S} \left \lbrack \sum_{s=1}^{S} \sum_{t\in \mathcal{T}} (1-m_{kts})X_{kts} \right \rbrack  \right) && \label{app-obj}\\
\text{s.t. }& 0 \leq V_{it} \leq \bar{v}_i(\kappa_t), \ &&  \forall i \in {\bf N}, t \in \mathcal{T},
\label{eqn:app-Speed_bound}\\
& 0 \leq I_{i0} \leq \bar{\iota}_i(\kappa_0),\ && \forall i \in {\bf E}^m, \label{eqn:app-Inventory-Bound-Init} \\
& X_{0ts} = \gamma_{0} d_{0ts}V_{1t}, \ &&\forall t \in \mathcal{T}, s \in \boldsymbol{\mathcal{S}} \label{eqn:app-System_feed}\\
& X_{its} = (1-\phi_i-\varphi_{it}(\kappa_t))\sum_{j \in \boldsymbol{\delta_i^-}}X_{jts} ,\ &&\forall i \in {\bf E}^g, t \in \mathcal{T}, s \in \boldsymbol{\mathcal{S}}  \label{eqn:app-Grinder_flow}\\
&X_{its} \leq  \gamma_{i}d_{its} V_{it}, \ &&\forall i \in {\bf E}^r,  t\in \mathcal{T}, s \in \boldsymbol{\mathcal{S}} \label{eqn:app-Conveyor_cap}\\
&X_{its} = \sum_{j \in \boldsymbol{\delta_i^-}}X_{jts},   \ &&\forall i \in {\bf E}^r \setminus \{b1, b2\},  t\in \mathcal{T}, s \in \boldsymbol{\mathcal{S}} \label{eqn:app-Conveyor_flow1} \\
&X_{b_1,ts} = (1 - \theta_{ts}) X_{bts},   \ &&\forall  t\in \mathcal{T}, s \in \boldsymbol{\mathcal{S}} \label{eqn:app-Conveyor_flow2}\\
&X_{b_2,ts} = \theta_{ts} X_{bts},  \ &&\forall  t\in \mathcal{T}, s \in \boldsymbol{\mathcal{S}} \label{eqn:app-Conveyor_flow3}\\
& X_{its} = \gamma_{i}d_{its}V_{it}, \ &&\forall i \in {\bf E}^m, t \in \mathcal{T}, s \in \boldsymbol{\mathcal{S}} \label{eqn:app-Metering_flow}\\
& I_{i1s} = I_{i0} + \sum_{j \in \boldsymbol{\delta_i^-}}X_{j1s} - X_{i1s}, \ &&\forall i \in {\bf E}^m, s \in \boldsymbol{\mathcal{S}}\label{eqn:app-MB_level1} \\
& I_{its} = I_{i(t-1)s} + \sum_{j \in \boldsymbol{\delta_i^-}}X_{jts} - X_{its}, \ &&\forall i \in {\bf E}^m, t \in \mathcal{T}\setminus\{1\}, s \in \boldsymbol{\mathcal{S}} \label{eqn:app-MB_level2} \\
& I_{its} \geq \underline{\iota}_i d_{its}, \ &&\forall i \in {\bf E}^m , t \in \mathcal{T}, s \in \boldsymbol{\mathcal{S}} \label{eqn:app-MB_cap1}\\
& I_{its}  \leq \bar{\iota}_i d_{its}, \ &&\forall i \in {\bf E}^m, t \in \mathcal{T}, s \in \boldsymbol{\mathcal{S}} \label{eqn:app-MB_cap2} \\
& X_{its} = (1-\varphi_{it}(\kappa_t)) \left( (1-\frac{\hat{t}_i}{t}) \sum_{j \in \boldsymbol{\delta_i^-}}X_{jts} + I_{i(t-1)s} \right), \ &&\forall i \in {\bf E}^{mill}, t \in \mathcal{T}, s \in \boldsymbol{\mathcal{S}}\label{eqn:app-Pellet_flow} \\
& I_{its} = I_{i(t-1)s} + \sum_{j \in \boldsymbol{\delta_i^-}}X_{jts} - X_{its}, \ &&\forall i \in {\bf E}^{mill}, t \in \mathcal{T}, s \in \boldsymbol{\mathcal{S}} \label{eqn:app-Pellet_level}\\
& I_{its}  \leq \bar{\iota}_i d_{its}, \ &&\forall i \in  {\bf E}^{mill}, t \in \mathcal{T}, s \in \boldsymbol{\mathcal{S}} \label{eqn:app-Pellet_cap} \\
& (1-m_{its})X_{its} \leq \bar{x}_i(\kappa_t), \ &&\forall i \in {\bf E}^p, t \in \mathcal{T},  s \in \boldsymbol{\mathcal{S}} \label{eqn:app-Eqn_Infeed_cap}\\
&(1- m_{kts}) X_{kts} \leq \overline{r}, \ &&\forall t \in \mathcal{T}, s \in \boldsymbol{\mathcal{S}} \label{eqn:app-Reactor_fed_cap}\\
& X_{0ts}, X_{its} \geq 0, \ &&\forall i \in {\bf N}, t \in \mathcal{T},  s \in \boldsymbol{\mathcal{S}}\\
& I_{i0}\geq 0, \ &&\forall i \in {\bf E}^m \cup {\bf E}^{mill},\\
& I_{its} \geq 0, \ &&\forall i \in {\bf E}^m \cup {\bf E}^{mill}, s \in \boldsymbol{\mathcal{S}}\\
\nonumber \\
&\displaystyle \frac{1}{S}\sum_{s=1}^S  \mathbb{1}\left[\frac{1}{T} \sum_{t \in \mathcal{T}} (1- m_{kts}) X_{kts}  \geq r \right]\geq 1-\hat{\epsilon}.\label{eqn:FeedTarget}
\end{alignat}
\end{footnotesize}
The bypass ratio in the separation unit is calculated as follows:
\begin{equation}\label{eq:bypassRatio}
\theta_{ts} := \begin{cases}
 \max\{0.5 - 0.4(\rho_{1ts}^{50}-6.35)/(\rho_{1ts}^{50} - \rho_{1ts}^{10}), 0\} &\text{for $\rho_{1ts}^{50} \geq 6.35$}\\
\min\{0.5 + 0.4(6.35 - \rho_{1ts}^{50})/(\rho_{1ts}^{90} - \rho_{1ts}^{50}), 1\} &\text{for $\rho_{1ts}^{50} < 6.35$},
\end{cases}
\end{equation}
where, the value 6.35 (in mm) corresponds to the screen size used in the separation unit. 

\pagebreak
We explain the constraints in formulation $(\hat{P})$ as follows:
\vspace{-2in}
\begin{itemize}
    \item Constraints \eqref{eqn:app-Speed_bound} and \eqref{eqn:app-Inventory-Bound-Init} correspond to the constraints \eqref{eqn:Speed_bound} and constraints \eqref{eqn:Inv_bound_init} in the succinct formulation ($\hat{P}$) in the main document. They represent the bounds on the equipment processing speed and initial inventory, respectively. 
    \item Constraints \eqref{eqn:app-System_feed}, \eqref{eqn:app-Grinder_flow}, \eqref{eqn:app-Metering_flow}, and \eqref{eqn:app-Pellet_flow} represent the flow calculations for processing and storage equipment. These constraints correspond to constraints \eqref{eqn:sFlow2} in the succinct formulation ($\hat{P}$) in the main document. For example, constraints \eqref{eqn:app-Grinder_flow} calculate the biomass flow with respect to the moisture and dry matter losses during the grinding process. 

    \item Constraints \eqref{eqn:app-Conveyor_cap} represent the upper limit on the amount of biomass flow from each transportation equipment $i \in {\bf E}^r$. These constraints correspond to constraints \eqref{eqn:sFlow3} in the succinct formulation ($\hat{P}$) in the main document.   

    \item In the transportation equipment, the biomass flowing  into the equipment equals to the biomass flowing from it. Constraints \eqref{eqn:app-Conveyor_flow1}, \eqref{eqn:app-Conveyor_flow2}, and \eqref{eqn:app-Conveyor_flow3} represent these flow balance equations, just like constraints \eqref{eqn:sFlow1} in the succinct formulation ($\hat{P}$) in the main document. 

    \item Constraints \eqref{eqn:app-MB_level1} and \eqref{eqn:app-MB_level2} are the inventory balance constraints of the storage equipment $i \in {\bf E}^m$. They are also a part of constraints \eqref{eqn:sFlow2} in the succinct formulation ($\hat{P}$) in the main document. 

    \item Constraints \eqref{eqn:app-MB_cap1} and \eqref{eqn:app-MB_cap2} set the upper and lower thresholds for the inventory level. The inventory upper bound comes from the storage equipment's capacity. 
    The lower threshold is required for consistent flow out of the storage equipment ${\bf E}^m$. Although the volume of a storage equipment is fixed, variations in the biomass density impact the amount of biomass allowed to be stored. These constraints are the detailed version of constraints \eqref{eqn:sInv_bound} in the succinct formulation ($\hat{P}$) in the main document.

    \item Pelleting process takes $\hat{t}_i$ units of time (which is less than the chosen time period $t$) in the pelleting mill $i$. The pelleting mill has an in-process storage capacity, which keeps the biomass during the pelleting process. Constraints \eqref{eqn:app-Pellet_level} and \eqref{eqn:app-Pellet_cap} represent the inventory balance and inventory capacity in pelleting mill ${\bf E}^{mill}$, respectively. In the succinct formulation in the main documents, we included these constraints within constraints \eqref{eqn:sFlow2} and \eqref{eqn:sInv_bound}, respectively.

    \item Processing equipment and the reactor have limits on the amount of biomass flowing into them. Constraints \eqref{eqn:app-Eqn_Infeed_cap} and \eqref{eqn:app-Reactor_fed_cap} represent these upper limits. These constraints are part of constraints \eqref{eqn:sFlow2} in the succinct formulation ($\hat{P}$) in the main document. 

\end{itemize}

\pagebreak
\section{Data Tables}\label{apdx:Cost}
\begin{table}[htp!]
	\scriptsize
	\centering
	{%
		\begin{tabular}{lcccccc}
			\toprule			
			& \multicolumn{2}{c}{\textbf{High moisture}}& \multicolumn{2}{c}{\textbf{Med moisture}}& \multicolumn{2}{c}{\textbf{Low moisture}}\\
			\cmidrule(lr){2-3}
			\cmidrule(lr){4-5}
			\cmidrule(lr){6-7}
			\textbf{Equipment} & \textbf{Energy} & \textbf{Fixed}& \textbf{Energy} & \textbf{Fixed}& \textbf{Energy} & \textbf{Fixed}\\
			& \textbf{Cost} & \textbf{Cost} & \textbf{Cost} & \textbf{Cost} & \textbf{Cost} & \textbf{Cost}\\ 
			& \textbf{$ (\$/hr)$}& \textbf{$ (\$/hr)$} & \textbf{$ (\$/hr)$} & \textbf{$ (\$/hr)$} & \textbf{$ (\$/hr)$} & \textbf{$ (\$/hr)$}\\
				\cmidrule(lr){1-1}
			\cmidrule(lr){2-2}
			\cmidrule(lr){3-3}
			\cmidrule(lr){4-4}
			\cmidrule(lr){5-5}
			\cmidrule(lr){6-6}
			\cmidrule(lr){7-7}
			Bale conveyor& 	$0.12	$& $0.48$&	$0.12	$&$0.48$&	$0.12	$&$0.48$\\
			Grinder 1 & $1.72	$ &$31.32$ &	$1.47$&	$31.32$	&$0.48$	&$31.32$\\
			Screw conveyor-6&	$0.29	$&$11.09$&	$0.29$	&$11.09$&	$0.29$&	$11.09$\\
			Drag chain conveyor-5&	$0.17	$&$1.18$&	$0.17$	&$1.18$	&$0.17$	&$1.18$\\
			Drag chain conveyor-6&	$0.29	$& $1.56$&	$0.29$ &	$1.56$ &	$0.29$ &	$1.56$\\
			Screw conveyor-1&	$0.17$ &	$1.20$&	$0.17$&	$1.20$ &	$0.17$ &	$1.20$\\
			Grinder 2&	$3.33$ &	$13.85$ &	$1.11$ &	$13.85$ &	$0.98$ &	$13.85$\\
			Screw conveyor-2&	$1.15$	&$4.22$	&$1.15$&	$4.22$&	$1.15$&	$4.22$\\
			Screw conveyor-4&	$0.29$&	$11.09$	& $0.29$ &	$11.09$ &	$0.29$	& $11.09$\\
			Metering bin&	$0.63$	&$9.98$	&$0.63$&	$9.98$	&$0.63$&	$9.98$\\
			Screw conveyor-5&	$0.29$&	$6.14$&	$0.29$&	$6.14$&	$0.29$	&$6.14$\\
			Pellet mill	& $6.06$&	$11.59$&	$4.07$&	$11.59$&	$3.70$&	$11.27$\\
			Drag chain conveyor-1& 	$0.12$	&$1.61$&	$0.12$&	$1.61$&	$0.12$&	$1.61$\\
			\midrule
			\textbf{Total}&  	$\boldsymbol{14.63}$	&$\boldsymbol{105.30}$&	$\boldsymbol{10.16}$&	$\boldsymbol{105.30}$&	$\boldsymbol{8.67}$&	$\boldsymbol{104.98}$\\			
			\bottomrule
	\end{tabular}}
	\caption{Energy Consumption and Fixed Costs}
	\label{tab:Cost}
\end{table}%

\begin{table}[htp!]
	\scriptsize
	\centering
	{%
		\begin{tabular}{lcccc}
			\toprule
			
			& \multicolumn{2}{c}{\textbf{Biomass Processed in Grinder 1}}& \multicolumn{2}{c}{\textbf{Biomass Processed in Grinder 2}}\\
			\cmidrule(lr){2-3}
			\cmidrule(lr){4-5}
			\textbf{Moisture} & $\boldsymbol{\rho^{50}}$& $\boldsymbol{\rho^{90}/\rho^{10}}$ & $\boldsymbol{\rho^{50}}$& $\boldsymbol{\rho^{90}/\rho^{10}}$\\
			\textbf{ \hspace{.1mm} Level}&$(mm)$ & &$(mm)$ &\\
			\cmidrule(lr){1-1}
			\cmidrule(lr){2-2}
			\cmidrule(lr){3-3}
			\cmidrule(lr){4-4}
			\cmidrule(lr){5-5}
			Low		&  $\left \lbrack 1.90 - 2.00\right \rbrack$ & $\left \lbrack 11.5 - 13.5 \right \rbrack$ & $\left \lbrack 0.70 - 0.60 \right \rbrack$& $\left \lbrack 7.00 - 6.00 \right \rbrack$ \\
			Medium	&  $\left \lbrack 2.30 - 2.40\right \rbrack$ & $\left \lbrack 11.0 - 13.0 \right \rbrack$ & $\left \lbrack 0.75 - 0.65 \right \rbrack$& $\left \lbrack 7.00 - 6.00 \right \rbrack$ \\
			High	&  $\left \lbrack 1.70 - 1.80\right \rbrack$ & $\left \lbrack 9.0 - 11.00 \right \rbrack$ & $\left \lbrack 0.65 - 0.55 \right \rbrack$& $\left \lbrack 9.00 - 8.00 \right \rbrack$ \\			
			\bottomrule
	\end{tabular}}
	\caption{Particle Size Distribution Percentiles}
	\label{tab:PSD}
\end{table}%

\begin{table}[htp!]
	\scriptsize
	\centering
	{%
		\begin{tabular}{lccccc}
			\toprule
			
			& \multicolumn{3}{c}{\textbf{Moisture Loss}}& \multicolumn{2}{c}{\textbf{Dry Matter Loss}}\\
			\cmidrule(lr){2-4}
			\cmidrule(lr){5-6}
			\textbf{Moisture} & \textbf{Grinder 1}& \textbf{Grinder 2} & \textbf{Pelleting Mill}& \textbf{Grinder 1} & \textbf{Grinder 2}\\
			\textbf{ \hspace{.1mm} Level}&$(\%)$ &$(\%)$ &$(\%)$ &$(\%)$ &$(\%)$\\
			\cmidrule(lr){1-1}
			\cmidrule(lr){2-2}
			\cmidrule(lr){3-3}
			\cmidrule(lr){4-4}
			\cmidrule(lr){5-5}
			\cmidrule(lr){6-6}
			Low		&  $0.50$ & $0.7 0$ & $0.00$ & $1.50$ & $0.50$ \\
			Medium	&  $3.00$ & $3.00$ & $1.50$ & $1.50$ & $0.50$ \\
			High	&  $4.77$ & $4.00$ & $3.90$ & $1.50$ & $0.50$ \\			
			\bottomrule
	\end{tabular}}
	\caption{Moisture and Dry Matter Changes}
	\label{tab:processing_changes}
\end{table}%

\begin{table}[htp!]
	\scriptsize
	\centering
	{%
		\begin{tabular}{lccc}
			\toprule			
			& \textbf{Low Moisture}& \textbf{Med Moisture}& \textbf{High Moisture}\\
			& $(dt/ hr)$& $(dt/ hr)$& $(dt/ hr)$\\
			\cmidrule(lr){2-2}
			\cmidrule(lr){3-3}
			\cmidrule(lr){4-4}
			Grinder 1 & 5.23&4.53 &2.20 \\
			Grinder 2 &5.23 &2.80 &1.59 \\
			Pelleting Mill &4.76 & 3.81& 3.34\\
			\bottomrule
	\end{tabular}}
	\caption{Equipment Infeed Rate Limits}
	\label{tab:infeed_limit}
\end{table}%

\pagebreak
\section{A Bisection Search Based Heuristic Algorithm for Solving the SAA Problem ($\hat{P}$)}\label{apdx:BinearySearch}
\begin{algorithm}[!htp]
\scriptsize
\caption{A bisection search based algorithm for solving the SAA problem ($\hat{P}$)}
{\bf Input:}  Set $\overline{\pi} = 10^7$, $\underline{\pi} = 0$, $\delta = 0.01$ and $\sigma = 0.01$ 
\label{Algo:binary}
\begin{algorithmic}[1]
\STATE Set $\pi \gets \overline{\pi}$.
\STATE Solve model ($\bar{P}$) and get an optimal solution $\hat{\mathcal{U}}_s$.
\STATE Set a counter $C \gets 0$.
\FOR {$s \in \mathcal{S}$}
\IF{$\hat{\mathcal{U}}_s > 0$}
\STATE $C \gets C_1+1$
\ENDIF
\ENDFOR
\IF{$C\geq\hat{\epsilon} N + \sigma$}
\STATE Status $\gets$ ``INFEASIBLE: Target R is not achievable".
\ELSE
\WHILE{$\frac{\overline{\pi}- \underline{\pi}}{2} > \delta$}
\STATE Set $\pi \gets \frac{\overline{\pi}+ \underline{\pi}}{2}$
\STATE Solve model ($\bar{P}$) and get an optimal solution $\hat{\mathcal{U}}_s$.
\STATE Set a counter $C \gets 0$.
\FOR {$s \in \mathcal{S}$}
\IF{$\hat{\mathcal{U}}_s > 0$}
\STATE $C \gets C+1$
\ENDIF
\ENDFOR
\IF{$C\geq\hat{\epsilon} N + \sigma$}
\STATE $\underline{\pi} \gets \frac{\overline{\pi}+ \underline{\pi}}{2}$
\ELSE
\STATE $\overline{\pi} \gets \frac{\overline{\pi}+ \underline{\pi}}{2}$
\ENDIF
\ENDWHILE
\STATE Let $\widehat{\mathbf{V}}$, $\widehat{\mathbf{X}}$, $\widehat{\mathbf{I}}$ be the optimal solution of model ($\bar{P}$) in the final iteration.
\STATE Return $\widehat{\mathbf{V}}$, $\widehat{\mathbf{X}}$, $\widehat{\mathbf{I}}$
\ENDIF
\end{algorithmic}\label{alg:BinearySearch}
\end{algorithm} 

\section{Regression Analysis for Biomass Density}
\begin{table}[htp!]
	\scriptsize
	\centering
	{%
		\begin{tabular}{lccccc}
			\toprule			
			& & \multicolumn{4}{c}{\textbf{P-value}}\\
			\cmidrule(lr){3-6}
			\textbf{Regression} & $\mathbf{R^2}$ & $\mathbf{\alpha_i^0}$ &$\mathbf{m_{it}}$ & $\mathbf{\rho^{50}_i}$ & $\mathbf{\frac{\rho^{90}_i}{\rho^{10}_i}}$\\
			\cmidrule(lr){1-1}
			\cmidrule(lr){2-2}
			\cmidrule(lr){3-3}
			\cmidrule(lr){4-4}
			\cmidrule(lr){5-5}
			\cmidrule(lr){6-6}
			(16) & 0.956 & $4.1*10^{-142}$ & $9.5*10^{-134}$ & $1.3*10^{-69}$ & $6.9*10^{-1}$\\
			(17) & 0.945 & $2.1*10^{-69}$ & $2.4*10^{-75}$ & $2.5*10^{-39}$ & $7.3*10^{-2}$\\
			\bottomrule
	\end{tabular}}
	\caption{Regression Analysis Statistics}
	\label{tab:reg_stats}
\end{table}%
\end{document}